\documentclass[letterpaper, paper,11pt]{AAS}		

\usepackage{bm}
\usepackage{amsmath}
\usepackage{caption}
\usepackage{subcaption}
\usepackage{mathtools}
\usepackage{float}

\usepackage[colorlinks=true, pdfstartview=FitV, linkcolor=black, citecolor= black, urlcolor= black]{hyperref}
\usepackage{overcite}
\usepackage{footnpag}			      	

\PaperNumber{23-313}

\begin{document}

\newpage
\thispagestyle{empty}
\mbox{} 

\noindent Disclaimer: 

This manuscript was presented at the 33rd AAS/AIAA Space Flight Mechanics Meeting, Austin, TX 2023 (paper ID:23-313). 
This version is now obsolete. Please reference the updated version in Celestial Mechanics and Dynamical Astronomy (Springer): \url{https://doi.org/10.1007/s10569-023-10179-8}.

\newpage

\title{Automated Tour Design in the Saturnian System}

\author{Yuji Takubo \thanks{Undergraduate Student, School of Aerospace Engineering, Georgia Institute of Technology, 30332.},  
Damon Landau \thanks{Mission Formulation Engineer, Systems Engineering Section, Jet Propulsion Laboratory, California Institute of Technology, 91011},
Brian Anderson \thanks{Mission Analyst, Mission Design and Navigation Section, Jet Propulsion Laboratory, California Institute of Technology, 91011},
}

\maketitle{} 		

\begin{abstract}
Future missions to Enceladus would benefit from multi-moon tours that leverage $V_\infty$ on resonant orbits to progressively transfer between moons. Such ``resonance family hopping" trajectories present a vast search space for global optimization due to the different combinations of available resonances and flyby speeds. The proposed multi-objective tour design algorithm optimizes entire moon tours from Titan to Enceladus via grid-based dynamic programming, in which the computation time is significantly reduced by utilizing a database of $V_\infty$-leveraging transfers. The result unveils a complete trade space of the moon tour design to Enceladus in a tractable computation time and global optimality.   
\end{abstract}

\section{Introduction}
Ocean worlds such as Europa, Titan, and Enceladus present ideal candidates to search for life outside of Earth. Such scientific interests are enabled by multi-moon flyby trajectories of planetary systems, which requires an optimization under a unique dynamical system\cite{mackenzie2021enceladus}. Here, the spacecraft (SC) performs multiple flybys (gravity assists, GA) to progressively reach inner moons via the so-called \textit{resonance family hopping}, or \textit{pump-down} trajectories. The designed trajectories are often plotted on Tisserand graphs (or Tisserand-Poincare leveraging graph, TG) \cite{strange2002tg,campagnola2010endgame1,campagnola2010endgame2}, which identify the properties of trajectories between flybys. 

Koon et al. optimized a Jovian multi-moon orbiter design based on the patched three-body approach  \cite{kolomaro2006ThreeBody}, while most of the previous literature had adopted methods that utilize $V_\infty$-leveraging transfers (VILT) \cite{Sims1997VILT}. Regarding the moon tour to Enceladus, Strange et al. \cite{strange2009enceladus}, Campangnola et al. \cite{campagnola2010enceladus}, Palma \cite{palma2016enceladus}, and Landau \cite{landau2018efficient} developed the state-of-the-art benchmarks, designing the flyby sequence at each major moon of Saturn: Titan, Rhea, Dione, Tethys, and Enceladus. 

Nonetheless, several key issues are yet to be addressed in the tour design problem.
First, established design procedures rely on heuristics of the mission designers in order to produce a single trajectory\cite{campagnola2019tour}. Dynamic programming (DP) techniques can globally search the state space, but have been limited to only a few ($\leq$ 3) flybys at a time \cite{brinckerhoff2009pathfinding}. 
Since some tours of Saturn's moons require dozens of flybys, a straightforward implementation of branch and bound suffers the curse of dimensionality. In this case, optimization can take multiple days of computation time per moon \cite{palma2016enceladus}. 
Landau \cite{landau2018efficient} alleviated this computational burden by introducing primer vector theory to each arc, while still taking 17 hours with a 6-core CPU to complete the broad search.   
Evolutionary algorithms can also perform a global optimization in a stochastic manner, with the state of the art producing tours with less than 10 flybys. \cite{Ellithy2022globalmoon}. 
Furthermore, each moon tour has been designed separately with fixed boundary conditions \cite{strange2009enceladus, campagnola2010enceladus, palma2016enceladus} or lingering discontinuities. \cite{brinckerhoff2009pathfinding}. In a realistic mission design, boundary conditions should also be treated as free variables for end-to-end optimization at the possible expense of additional computation time. 

We address successive multi-moon tour design based on two key approaches. 
First, the \textit{path-finding} algorithm optimizes the sequences of resonance trajectories through grid-based DP \cite{bellome2022dp}. The deterministic optimization schemes have been adopted by previous literature \cite{landau2019global, Kaela2019europa}. A new graphical analysis tool, the pump-$V_\infty$ map, displays each resonance family as a curve of ballistic resonant orbits. The discretized points on these curves are then connected on the basis of the branch and bound method, which recursively propagates the resonance family hopping tour in the manner of DP (see references \cite{Beckman1974pump,Uphoff1976pump} or Figure \ref{fig:leg_opt_nomen} for definition of pump angle and $V_\infty$).
Second, for the leg optimization problem (i.e., \textit{path-solving}), we generate a database that covers all potentially useful combinations of resonance and $V_\infty$. Each same-moon leg is approximated to a linear model and captures the sensitivity of $\Delta V$ to optimally change $V_\infty$. This database is interpolated during path-finding to provide the optimal VILT from the given initial and terminal $V_\infty$, which allows us to approximate optimal legs without solving the computationally expensive optimal control problem during the path-finding. 

The contribution of this paper is threefold. First, we automate the solution of the entire multi-moon tour problem to be completed in only less than 30 \textcolor{blue}{minutes} with minimal interaction. This reduction in run time is made possible by discretization of both path-finding and path-solving variables to maintain a tractable search space on each successive flyby\textcolor{blue}{, as well as the precomputation of the sets of feasible transfers from each (discretized) resonance orbit to the other}. 
We also validate that the obtained solution sets are accurate enough to be further optimized by a high-fidelity trajectory design tool.
Second, the pump-$V_\infty$ map provides new insights into the general problem of transfer between resonant orbits. 
We show that the coordinate system using the pump angle and $V_\infty$ has an elegant connection to the representation of the underlying optimal control problem. 
Finally, our deterministic multi-objective optimization reveals an unexplored trade space of the end-to-end Saturn-moon tours \textcolor{blue}{with guaranteed optimality.}

\section{Problem formulation}

We consider the Saturn tours that begin at Titan and successively fly by Rhea, Dione, and Tethys before arriving at Enceladus. Table \ref{tab:moon_param} summarizes the properties of these moons and indicates the minimum allowable flyby altitudes. 
Since the eccentricity and inclination of each moon is small, we assume the moons follow concentric and co-planar circular orbits.
Similar to other multiple gravity assists (MGA) problems, the moon tour optimization problem has a hierarchical structure that is composed of two main parts: path-finding and leg optimization. 

\begin{table}[h!]
\centering
\begin{tabular}{c c c c c c} 
 \hline
 & Titan & Rhea & Dione & Tethys & Enceladus \\ 
 \hline\hline
 $a$ (km) & 1221870 & 527108 & 377396 & 294619 & 237948\\ 
 $e$ &  0.0288 & 0.001 & 0.0022 &  0.0001 & 0.0047 \\
 $i$ ($^{\circ}$) &  0.33 & 0.35 & 0.02 & 1.09 &  0.02  \\
 $r_M$ (km) & 2574.7 & 763.8 & 561.4 & 531.1 & 252.1\\
 period (days) & 15.945 & 4.152 & 2.737 & 1.89 & 1.370\\
 $GM$ $(\text{km}^3/\text{s}^2$) & 8977.9  & 153.94 & 73.110 & 41.209 & 7.2094 \\
 min. flyby altitude (km) & 1600 & 50 & 50 & 50 & 25 \\
 \hline
\end{tabular}
\caption{Parameters of Saturn's moon}
\label{tab:moon_param}
\end{table}

\begin{figure}[ht]
     \centering
     \begin{subfigure}[b]{0.37\textwidth}
         \centering
         \includegraphics[width=\textwidth]{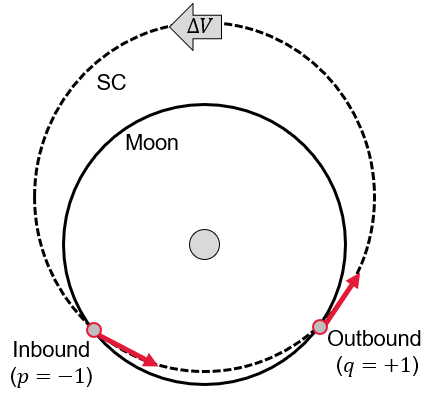}
         \caption{Exterior transfer $(M>N)$}
         \label{fig:exterior}
     \end{subfigure}
     \begin{subfigure}[b]{0.37\textwidth}
         \centering
         \includegraphics[width=\textwidth]{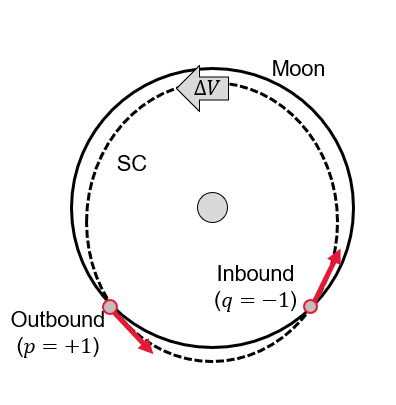}
         \caption{Interior transfer $(M<N)$}
         \label{fig:interior}
     \end{subfigure}
     \caption{Nomenclature in exterior and interior transfers}
\label{fig:exter_inter}
\end{figure}

First, the path-finding problem selects the optimal sequence of resonant transfers, which are elliptic orbits that re-encounter the same moon after several revolutions. Such orbits are expressed with four parameters $[M,N,p,q]$; we call this set of four parameters the \textit{resonance family}. The first two numbers indicate the (positive integer) resonance ratio of the moon revolution to the SC revolution, $M:N$. If $M > N$, the period of the SC is longer than that of the moon, i.e., the semimajor axis of the SC is larger than that of the moon, and vice versa. $M > N$ transfers are therefore called \textit{exterior} transfers, and $M < N$ transfers are called \textit{interior} transfers, as shown in Figure \ref{fig:exter_inter}. Transfers such that $M=N$ could be either interior or exterior transfers. 
It is optimal to perform the leveraging impulse near apoapsis for exterior transfers and near periapsis for interior transfers \cite{landau2018efficient}. 
The $p,q$ parameters indicate whether the SC encounters the moon inbound (-1) or outbound (+1) at the starting ($p$) and ending ($q$) points. Note that due to their symmetry, $[M,N,+1,+1]$ and $[M,N,-1,-1]$ ballistic trajectories have the exact same time of flight (ToF) and the absolute value of the pump angle. We retain the nomenclature \textit{resonance family} when $p \neq q$ even though those transfers are not exactly resonant (non-integer number of revolutions). 

\begin{figure}[ht]
  \centering
  \includegraphics[width=0.85\linewidth]{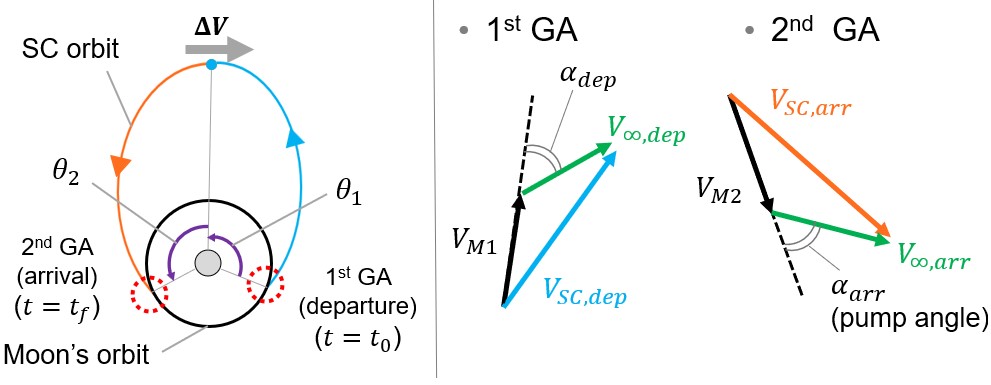}
  \caption{Nomenclature of resonant orbits for leg optimization. For simplicity, the shown SC trajectory has one revolution. In this example, $p = +1$ and $q = -1$.}
  \label{fig:leg_opt_nomen}
\end{figure}

Second, given the parameters of the resonance family in the path-finding problem, we must optimize each arc of the tour to traverse resonances with minimum $\Delta V$; we call this process as leg optimization. We allow one impulsive maneuver per leg, similar to the MGA-1DSM model \cite{vasile2006mga,englander2012auto}, as shown in Figure \ref{fig:leg_opt_nomen}. By solving the multi-revolution two-point boundary value problem (TPBVP) for each transfer, we propagate the resonant orbits, formulating the entire moon tour. 

While most trajectory optimization problems are solved as a single-objective problem, mission designers may want to see a multi-objective trade space, especially at an early stage of a mission design. In this problem, the minimization of total ToF and total $\Delta V$ are chosen as the objectives. 


\section{Path-finding problem: Grid-based Dynamic Programming}
\textcolor{blue}{A ballistic transfer has one more degree of freedom besides $[M,N,p,q]$, in which we use the the encounter speed $V_{\infty,arr}$  in this paper to uniquely define the transfer.} is also a free variable. An issue we now have is that $[M,N,p,q, V_{\infty,arr}]$ is no longer a set of discrete variables because $V_{\infty,arr}$ is continuous. To keep the search space tractable, we discretize the range of feasible $V_{\infty,arr}$ for each resonance family. This process provides a finite number of grids in the entire design space, which are then optimized by DP. We tune the grid size so that the solution remains accurate for the early-stage design, while reducing the computational cost compared to methods that maintain continuous parameters. 

\subsection{pump-$V_{\infty}$ map}
\begin{figure}[ht]
  \centering
  \includegraphics[width=0.9\linewidth]{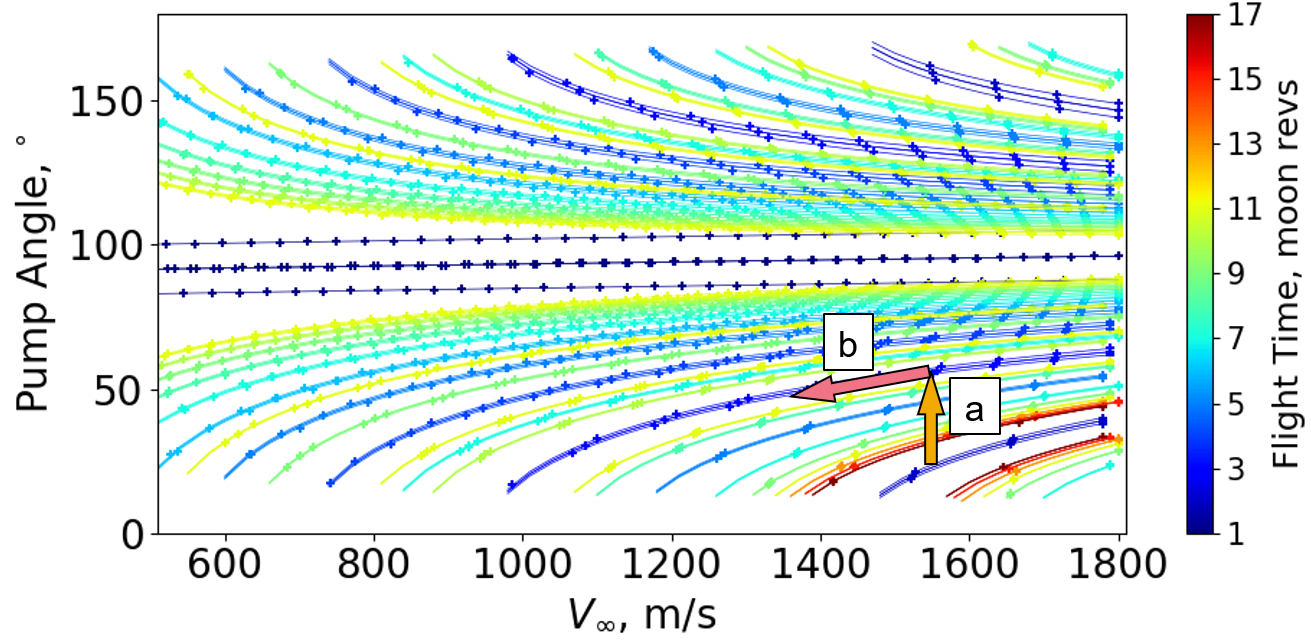}
  \caption{pump-$V_\infty$ map of Rhea $(1\leq M \leq17)$ with curves of $[M,N,p,q, V_{\infty,arr}]$}
  \label{fig:pump_vinf}
\end{figure}

We propose a new graphing tool, the \textit{pump-$V_\infty$ map}, which enables us to project the $[M$,$N$,$p$,$q$, $V_{\infty,arr}]$ parameters onto 2D space. A sample pump-$V_\infty$ map of Rhea is shown in Figure \ref{fig:pump_vinf} where the absolute values of the pump angles are shown on the map for convenience. 
Each point on this map represents a ballistic transfer, where each resonance family creates a curve colored by flight time. 
It is also observed that the same resonance $M:N$ with different $p$ and $q$ usually overlaps; as $M$ and $N$ become smaller, the discrepancy of the four curves becomes more pronounced, best represented by the $1:1$ resonance families which lie around $\alpha = 85 ^{\circ}, 90^{\circ},$ and $105^{\circ}$. 
The SC period increases with decreasing pump angle or increasing $V_\infty$, and decreases with increasing pump angle or decreasing $V_\infty$. 
Most commonly, a moon tour goes from the bottom right (high $V_\infty$ and small pump angle) to the top left (low $V_\infty$ with pump angle close to $180^{\circ}$) to traverse from the outer to inner moons. 
The tick marks along the resonance curves demonstrate the change in $V_\infty$ available from a constant $\Delta V$ (15 m/s in Figure \ref{fig:pump_vinf}), where resonances close to 1:1 are the least effective, and the $\partial V_\infty/\partial \Delta V$ efficiency improves with longer or shorter periods of the SC (i.e., bottom-right and top-right of the map, respectively).

\begin{figure}[ht]
     \centering
     \begin{subfigure}[b]{0.51\textwidth}
         \centering
         \includegraphics[width=\textwidth]{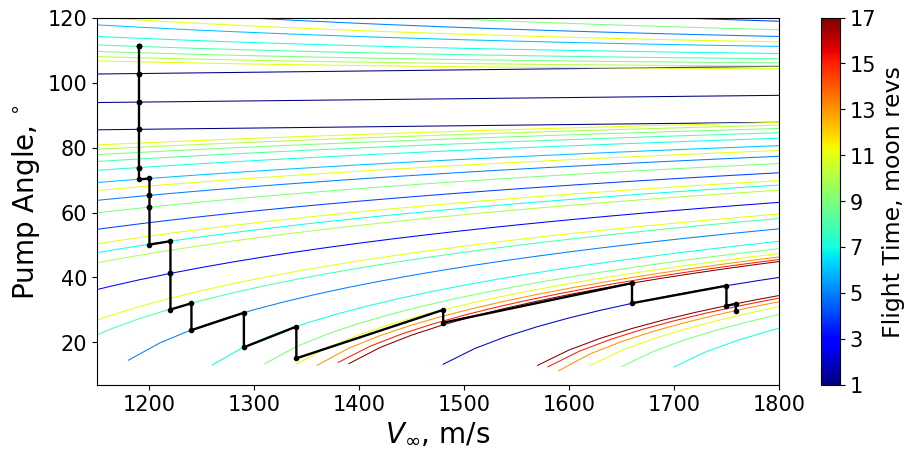}
         \caption{pump-$V_\infty$ representation}
         \label{fig:sample_pump_vinf}
     \end{subfigure}
     \begin{subfigure}[b]{0.46\textwidth}
         \centering
         \includegraphics[width=\textwidth]{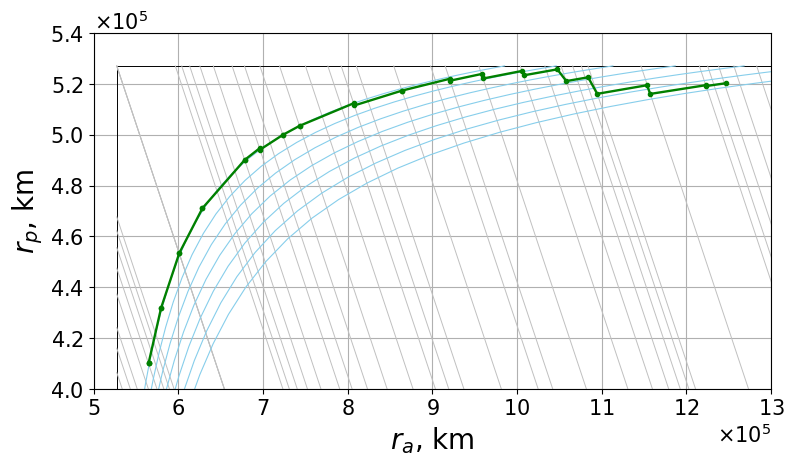}
         \caption{$r_p$ - $r_a$ (Tisserand graph) representation}
         \label{fig:sample_tg}
     \end{subfigure}
     \caption{Example Rhea tour design and comparison of the representation methods}
\label{fig:pump_vinf_and_tg}
\end{figure}

Any transfer can be decomposed into a combination of two maneuvers on this map. First, a vertical arrow (a) in the map indicates a (ballistic) flyby that keeps $V_\infty$ constant while changing the pump angle. This maneuver is responsible for switching from one resonance family to the next. Secondly, arrow (b) represents a change in $V_\infty$ along a resonance family curve, i.e., a VILT. 
A TG can also be used to represent the moon tour in 2D space. A sample Rhea tour plotted in the pump-$V_\infty$ map and $r_p-r_a$ TG are compared in Figure \ref{fig:pump_vinf_and_tg}. In Figure \ref{fig:sample_tg}, grey lines indicate resonance families, and blue curves indicate contours of the constant $V_\infty$ from 1800 m/s (right bottom) to 1200 m/s (left top) with interval of 100 m/s. 

The pump-$V_\infty$ map is tailored to our specific problem of tour design with a single moon and is not intended for tours that often switch between multiple moons. Philosophically, the pump-$V_\infty$ map focuses the design on flyby coordinates (in essence, a 2D representation of a $V_\infty$ globe\cite{strange2008globe}), whereas the $r_p-r_a$ TG emphasizes orbital coordinates. In both, the basic process is to sequentially add a flyby followed by an orbital transfer, with changes in $V_\infty$ confined along resonance family contours to maintain phasing. In pump-$V_\infty$ space the flyby occurs along a natural coordinate, where limits on maximum $\Delta \alpha$ per flyby can be represented by additional contours (or by directly scaling the pump axis) of number of flybys. In $r_p-r_a$ TG space flybys follow contours of constant $V_\infty$ superimposed on the map, where limits on maximum $\Delta \alpha$ per flyby can be represented by ticks spaced along each contour. In both, the resonance families are tracked as additional contours, where the $\Delta V$ cost to change $V_\infty$ can be represented by ticks spaced along each contour. The resolution of the flyby $V_\infty$ is arbitrary in pump-$V_\infty$ space, whereas many additional contours become necessary to track fine changes in $r_p-r_a$ space. We find that comparing both the flyby-centered and orbital-centered views of each tour provides a more complete understanding of the underlying mechanics.  

\subsection{Grid-based DP}
\begin{figure}[ht]
    \centering
    \begin{subfigure}[b]{0.48\textwidth}
        \centering
        \includegraphics[width=\textwidth]{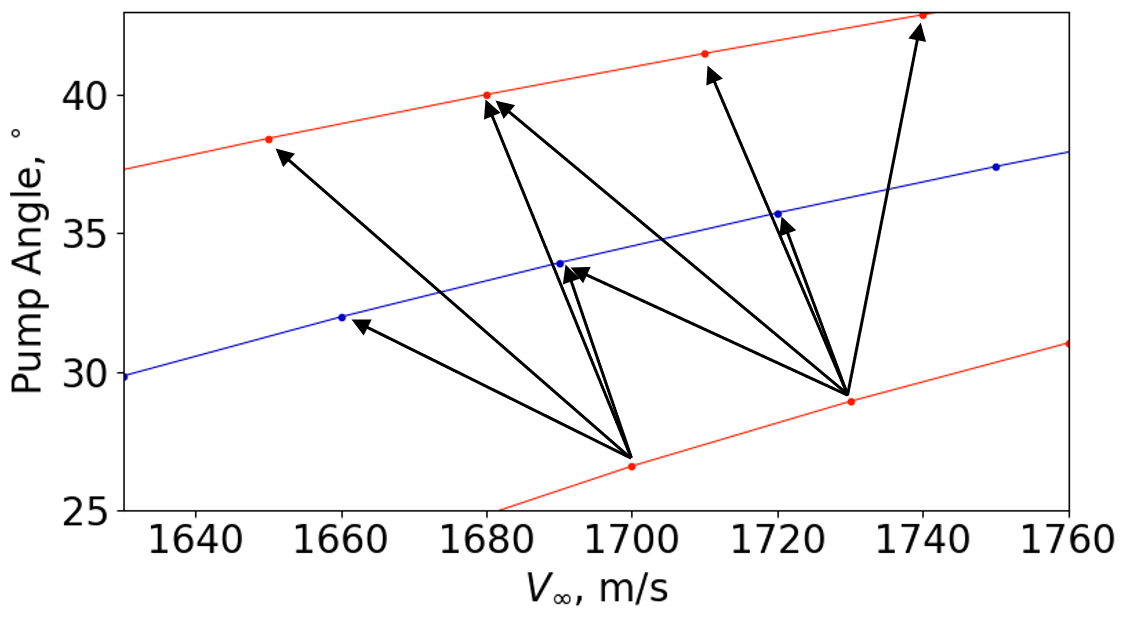}
        \caption{Branching out to the all feasible nodes}
        \label{fig:branch}
    \end{subfigure}
    \begin{subfigure}[b]{0.48\textwidth}
        \centering
        \includegraphics[width=\textwidth]{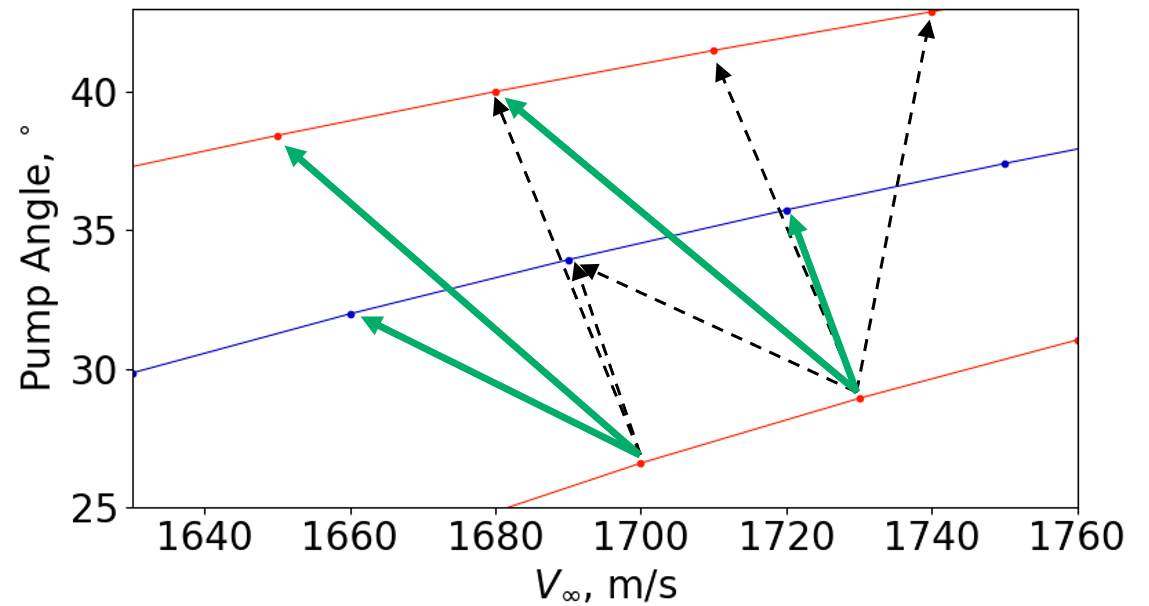}
        \caption{Pruning the solutions by 2D Pareto sorting}
        \label{fig:bound}
    \end{subfigure}
    \caption{Branch and bound in the grid-based DP. Black solid arrows, black dotted arrows, and green arrows represent reachable transfers, inferior (pruned) transfers, and Pareto-optimal transfers.}
\label{fig:b_and_b}
\end{figure}

The procedure of DP to solve the multi-objective path-finding problem is illustrated in Figure \ref{fig:b_and_b}. 
First, for each departing node from the initial node set, all other nodes on the pump-$V_\infty$ map are swept to check whether they can be reached from it. This way, we can branch out to all candidates of the \textcolor{blue}{arriving node for each departing node} (Figure \ref{fig:branch}). 
Note that the discretization of $V_{\infty, arr}$ enables us to limit the combination of \textcolor{blue}{departing and arriving} nodes. 
Then, we sort the feasible transfers based on the Pareto optimality, and prune all inferior solutions at each departing node after each maneuver (Figure \ref{fig:bound}).
We minimize two competing objectives: total ToF and total $\Delta V$ at each departing node.


After pruning \textcolor{blue}{the inferior arriving nodes at each departing node}, the set of Pareto-optimal arriving nodes is stored for use as the \textcolor{blue}{departing} nodes of the next maneuver. 
By recursively performing these branching and pruning processes, \textcolor{blue}{we obtain a large bundle of trajectories which connect locally Pareto-optimal transfers for each departing node. At the end of the whole moon tour(s), we can perform the Pareto-sorting for the all feasible trajectories but not each final node, which provides the set of globally Pareto-optimal moon tours.}  
Note that this optimization process is deterministic, so we obtain the exact Pareto front with guaranteed reproducibility.

\section{Leg Optimization: Linearized optimal VILT model}

\subsection{Multi-revolution TPBVP}
The leg optimization is formulated as a multi-revolution TPBVP. The problem is transcribed in Eq. \ref{eq:tpvbp}.

\begin{equation} \label{eq:tpvbp}
\begin{split}
 \text{given} & \quad t_0, V_\infty, \Delta V_\infty\\
 \min & \quad |\Delta V \left( \theta_1, \theta_2, \alpha_{dep}, \alpha_{arr}, t_f  \right)| = |\mathbf{v}(t_f - \Delta t_2) - \mathbf{v}(t_0 + \Delta t_1)| \\
 \text{subject to} & \quad \mathbf{r}(t_0 + \Delta t_1) = \mathbf{r}(t_f - \Delta t_2)\\
 & \quad t_0 + \Delta t_1 = t_f - \Delta t_2 \\
 & \quad \mathbf{r}(t_0) = \mathbf{r}_{M}(t_0) \\
 & \quad \mathbf{r}(t_f) = \mathbf{r}_{M}(t_f) \\
\end{split}
\end{equation}

\noindent where $V_{\infty}$ and $\Delta V_{\infty}$ are related to $V_{\infty, dep}$ and $\Delta V_{\infty, arr}$ as
\begin{equation} \label{eq:DVinf}
\begin{split}
    V_{\infty,dep} &= V_\infty + \Delta V_\infty \\ 
    V_{\infty,arr} &= V_\infty - \Delta V_\infty 
\end{split}
\end{equation}

\noindent Note that $V_{\infty,dep} = V_{\infty,arr}$, $\alpha_{\infty,dep} = \alpha_{\infty,arr}$, and $\Delta V_\infty = 0$ hold for ballistic transfers. Also, we can set $t_0 = 0$ without loss of generality. 

\begin{figure}[ht]
  \centering
  \includegraphics[width=0.8\linewidth]{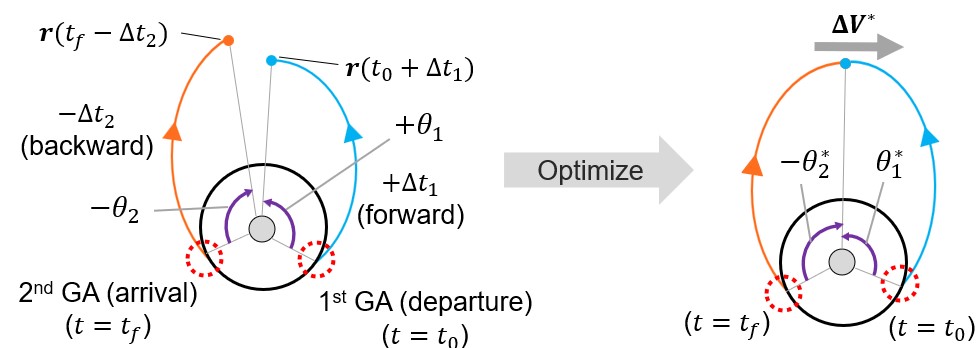}
  \caption{Nomenclature of the TPBVP.}
  \label{fig:tpbvp}
\end{figure}

The optimization scheme is shown in Figure \ref{fig:tpbvp}. First, we propagate the state from $\mathbf{r}(t_0)$  through a transfer angle of $\theta_1$, and propagate the state backwards from $\mathbf{r}(t_f)$ through an angle of $-\theta_2$. After these propagations, there would be a residual error in terminal positions and times that must be corrected, while the difference between the two velocity vectors creates an impulsive $\Delta V$ to be minimized. This problem is solved via the interior point method and accounts for $J_2$ perturbation from Saturn\cite{danielson1990j2}.

Since this is a multi-modal problem, an initial guess of the variables $\left( \theta_1, \theta_2, \alpha_{dep}, \alpha_{arr}, t_f  \right)$ critically affect the solution. Given the parameters $[M,N,p,q, V_{\infty,arr}]$, we use the ballistic solution as an initial guesses for $\alpha$ and $t_f$. The initial guess to place the maneuver ($\theta_{1,2}$) on VILTs is at apoapsis for exterior transfers and $M=N, (p,q)=(-1,+1)$, or at periapsis for interior transfers and $M=N, (p,q)\neq(-1,+1)$ (see Appendix, Eq. \ref{eq:f_to_thetas}). For $N>1$ we heuristically place the maneuver on the middle revolution, noting that the efficiency $\partial V_{\infty}/\partial \Delta V$ is the same on each rev in the linear model \cite{landau2018efficient}.

\subsection{Database approach and linearized trajectory model}

We can optimize a VILT for any combination of $[M,N,p,q]$ by solving Eq. \ref{eq:tpvbp}. 
However, solving this optimization problem for each leg at each iteration would be computationally expensive. 
Instead, we precompute all transfers likely to produce an optimal tour according to the design bounds specified in Table \ref{tab:vinf_ub}.
The key philosophy is to develop a database of such optimal VILTs and a corresponding linearized trajectory model by interpolating the database. 

The maximum number of moon revolutions $M$ sets an upper limit on flight time, and the $V_\infty$ bounds implicitly constrains the minimum and maximum number of SC revolutions $N$ and SC period. The Hohmann transfer $V_\infty$ and periods of the previous and following moons roughly bound the design space. 
The generation of the VILT database for all five moons completes in less than 15 minutes using an i7-8550U @ 1.80GHz (4-core) CPU with 16 GB RAM.

\begin{table}[h!]
\centering
\begin{tabular}{c c c c c c} 
\hline
 & Titan & Rhea & Dione & Tethys & Enceladus \\ 
\hline\hline
min $V_\infty \text{(m/s)} $  & 1200 & 650 & 550 & 550 & 200 \\ 
max $V_\infty \text{(m/s)} $  & 1600 & 1900 & 1000 & 900 & 850 \\ 
max $M$ & 2 & 15 & 15 & 16 & 25 \\ 
\hline
\end{tabular}
\caption{Bounds on $V_\infty$ and $M$ define the domain of the search space.}
\label{tab:vinf_ub}
\end{table}

Each row of the database contains three categories of information: 
\begin{itemize}
    \item Input parameters: $[M,N,p,q,V_\infty]$
    \item ToF and pump angle of ballistic transfer: $[t_f, \alpha]$
    \item First-order partial derivatives of properties of the leg with respect to optimal $\Delta V$: \\ $\frac{\partial}{\partial \Delta V} [t_f, V_{\infty,dep}, V_{\infty, arr}]$
\end{itemize}
To obtain these values, we perform two optimization problems for each input; $\Delta V_\infty = 0$ m/s is solved to obtain the ballistic transfer, and the partial derivatives are numerically approximated by solving another problem with small nonzero $\Delta V_\infty$ (5 m/s in this study) then taking the difference. 

\begin{figure}[ht]
     \centering
     \begin{subfigure}[b]{0.48\textwidth}
         \centering
         \includegraphics[width=\textwidth]{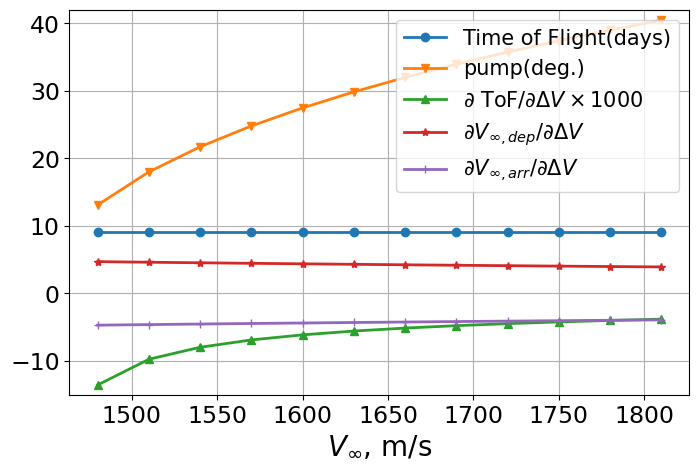}
         \caption{Graphs of the values of VILT database}
         \label{fig:graph1}
     \end{subfigure}
     \begin{subfigure}[b]{0.48\textwidth}
         \centering
         \includegraphics[width=\textwidth]{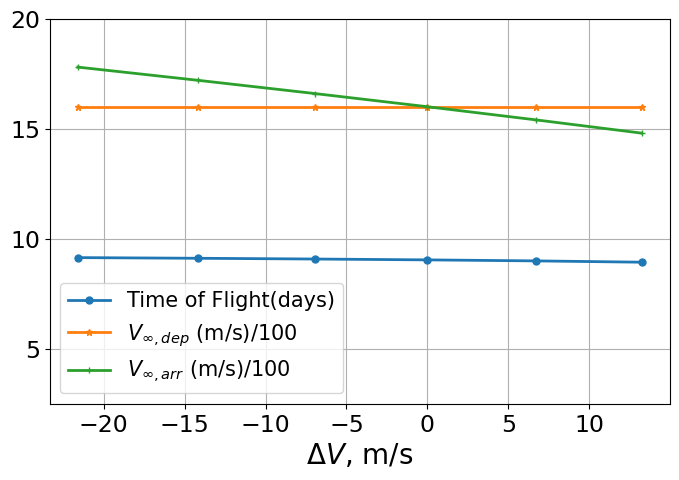}
         \caption{Reformulated graph with $V_{\infty, dep}=1600$ m/s}
         \label{fig:graph2}
     \end{subfigure}
     \caption{Sample Graphs of the trajectory: Rhea [2,1,+1,+1] family}
\label{fig:graph1_graph2}
\end{figure}

The partial derivatives are used to approximate all non-ballistic transfers up to $\pm$ 50 m/s $\Delta V$.  For each resonance family, we create piecewise continuous polynomials of the database as a function of $V_{\infty}$. Figure \ref{fig:graph1} shows a sample set of graphs for the Rhea  [2,1,1,1] family. Since $V_{\infty,dep}$ must be equal to $V_{\infty,arr}$ of the previous leg for a ballistic flyby, we can compute $\Delta V$ needed for each grid point of $V_\infty$ in Figure \ref{fig:graph1} to satisfy such a $V_{\infty,dep}$ as follows. 

\begin{equation} \label{eq:compute_dv}
\begin{split}
    \Delta V = \frac{V_{\infty, dep} - V_\infty}{\frac{\partial V_{\infty,dep}}{\partial \Delta V}}
\end{split}
\end{equation}

\noindent Note that by definition, positive $\Delta V$ increases $V_{\infty,dep}$ and decreases $V_{\infty, arr}$, and vice versa. 

The obtained $\Delta V$ is then used with other partial derivatives to compute the ToF and $V_{\infty,arr}$ as shown in Figure \ref{fig:graph2}. After computing the new graphs, we make grids of $V_{\infty,arr}$ corresponding to the nodes of the grid-based DP. 

Multiple constraints are enforced during the branching process in the grid-based DP, thus not all resonance orbits shown in Figure \ref{fig:graph2} can be reached from the current node. First, the range of reachable pump angle is limited by the maximum bend angle $|\Delta \alpha|$ at each flyby.
\begin{equation} \label{eq:max_pump}
\begin{split}
    |\Delta \alpha| \leq 2 \sin ^{-1} \frac{1}{1+\mu V_\infty ^2}
\end{split}
\end{equation}

\noindent where parameter $\mu$ is defined as
\begin{equation} \label{eq:mu}
\begin{split}
    \mu = \frac{r_M + \text{(min. flyby altitude)}}{GM_M}.
\end{split}
\end{equation}

\noindent This sets the upper and lower bounds of $\alpha_{dep}$ as follows.
\begin{equation} \label{eq:pump_boundary}
\begin{cases}
    \max \left( 0^{\circ},\alpha_{arr}^k - |\Delta \alpha|\right) \leq \alpha_{dep}^{k+1} \leq \min \left(180^{\circ}, \alpha_{arr}^k + |\Delta \alpha| \right) & \text{if} \quad \alpha_{arr}^k \geq 0 \\
    \max \left( -180^{\circ},\alpha_{arr}^k - |\Delta \alpha|\right) \leq \alpha_{dep}^{k+1} \leq \min \left(0^{\circ}, \alpha_{arr}^k + |\Delta \alpha| \right) & \text{if} \quad \alpha_{arr}^k < 0 \\    
\end{cases}
\end{equation}

\noindent where $k$ indicates the $k$-th transfer. Here we constrain the pump angle in order to prevent the SC from switching inbound/outbound direction (i.e., $p$ and $q$). Furthermore, we set a constraint on the $V_{\infty,arr}$ such that it does not exceed the lower bound of the $V_\infty$ of the given resonance family ($\alpha = 0^\circ$ for exterior transfers and $\alpha = 180^\circ$ for interior) to maintain the fidelity of the linearized trajectory model. 

\textcolor{blue}{Note that $\alpha_{arr}$ and $\alpha_{dep}$ would also be influenced by adding $\Delta V$. However, in this paper we assume that the pump angles are the function of solely the parameters $[M, N, p, q, V_\infty]$, but not $\Delta V$ for the following two reasons. 
First, this assumption made us to uniquely identify the orbit properties in the pump-$V_\infty$ map (i.e., all orbits after VILTs exactly locates on the curves of the ballistic transfers), which enables us to precompute the all combinations of the feasible transfers before the grid-based DP. 
We find that this precomputation significantly reduces the computation time of the DP, avoiding to compute the possible transfer at each manauver. 
Furthermore, this reduction in the trajectory model is still accurate enough for the preliminary design, and we confirmed that the optimized trajectory serves as a good enough initial guess for the high-fidelity solver, as shown in the later study.}

\section{Heuristics for the optimization}

The path-finding algorithm with a precomputed VILT database theoretically finds the Pareto optimal solution set. However, several heuristics are introduced to further reduce the problem size. 

First, we only consider $p = q = +1$ except for 1:1 resonance. As shown in Figure \ref{fig:pump_vinf}, the four curves of each $M:N$ resonance family mostly overlap when taking the absolute values of the pump angle. Thus, we can reduce the search space significantly by only considering one of them, and patch the boundary condition with $|\alpha_{arr}^k|$ and $|\alpha_{dep}^k+1|$; Eq. \ref{eq:pump_boundary} can be rewritten as follows. 
\begin{equation} \label{eq:pump_boundary_abs}
\begin{split}
    \max \left( 0^{\circ},|\alpha_{arr}^k| - |\Delta \alpha|\right) \leq |\alpha_{dep}^{k+1}| \leq \min \left(180^{\circ}, |\alpha_{arr}^k| + |\Delta \alpha| \right)
\end{split}
\end{equation}

\noindent We reconstruct the tour from the optimization result by choosing the appropriate sign of the pump angles so that they match $p$ and $q$ at departure and arrival, respectively. 
For 1:1 resonant orbits, we consider three orbits $[M,N,p,q]=$ [1,1,+1,+1], [1,1,+1,-1], and [1,1,-1,+1] as these three are distinctively different, while [1,1,+1,+1] and [1,1,-1,-1] completely overlaps near $\alpha = 90^{\circ}$ due to symmetry. 

\textcolor{blue}{We reduce the number of nodes to be visited at each flyby by binning the values of two objectives used for the 2D Pareto sorting at each departing node. At each departing node, many solutions nearly overlap when plotting the Pareto-optimal solutions after each transfer. 
Motivated by this observation, we make a grid space of this 2D space so that we can maintain diversity of the Pareto solution set, while reducing the number of arriving nodes to be propagated as the departing points at the next round of transfers. In this study, the intervals of the grid bins are defined as follows: $\delta \Delta V = 0.1$ m/s, and $\delta \text{ToF} = 2$ days.}

Finally, we set an upper bound to the magnitude of an impulsive maneuver s.t. \textcolor{blue}{ $-20 \leq \Delta V \leq 50$ m/s}. This is based on the empirical observation that a high $\Delta V$ usually does not produce an (Pareto-) optimal solution in this particular problem.

\subsection{Terminal condition of the path-finding: Titan - Tethys tour}
\begin{figure}[ht]
  \centering
  \includegraphics[width=0.8\linewidth]{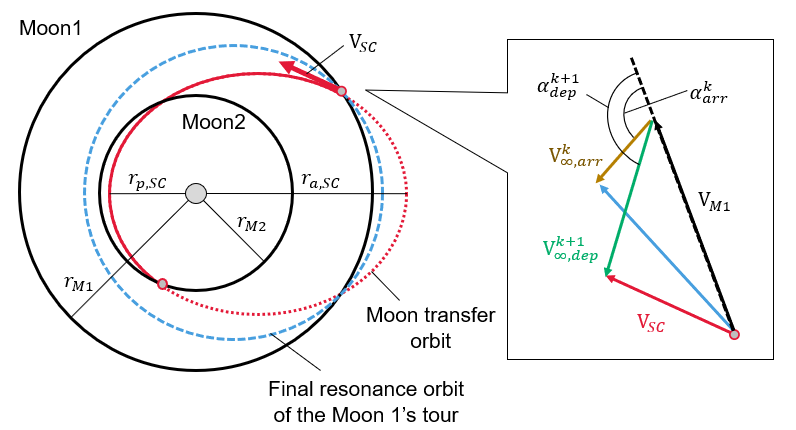}
  \caption{Transfer from a moon tour to the other}
  \label{fig:exit_nomen}
\end{figure}

In the previous literature \cite{strange2009enceladus, campagnola2010enceladus}, both the initial and terminal conditions \textcolor{blue}{(i.e., $V_\infty$ and $\alpha$)} of each moon tour were fixed. Instead of such a ``hard" criteria, we apply a ``soft" exit condition so that we can carry all the paths that can proceed to the next moon tour. A sample exiting transfer is shown in Figure \ref{fig:exit_nomen}.

Suppose that after the $k$-th flyby, a SC has an arriving pump angle $\alpha_{arr}^{k}$ and excess velocity $V_{\infty,arr}^{k}$. We assume that an unpowered $k+1$-th flyby with the maximum bend angle (Eq. \ref{eq:max_pump}) is performed (i.e., the minimum flyby altitude), which provides $\alpha_{dep}^{k+1}$. Using the law of cosine, the SC velocity is given as follows.

\begin{equation} \label{eq:V_sc}
    V_{SC} = \sqrt{V_{M1} ^ 2 + {V_{\infty,arr}^{k}} ^2 + 2  V_{M1} V_{\infty,arr}^{k} \cos \alpha_{dep}^{k+1}}
\end{equation}

\noindent where $V_{M1}$ is the velocity of the first moon and $V_{\infty,arr}^{k} = V_{\infty,dep}^{k+1}$. Also, we can obtain the angular momentum of the SC after the flyby as follows. 
\begin{equation} \label{eq:h_sc}
\begin{split}
    h_{SC} &= r_{M1} V_{SC,T} \\
           &= r_{M1} \left( V_{M1} + V_{\infty,arr}^{k} \cos \alpha_{dep}^{k+1} \right)
\end{split}
\end{equation}

\noindent where $V_{SC,T}$ indicates the transverse component of $V_{SC}$. Therefore, the semi-major axis and eccentricity of the SC orbit is expressed as 
\begin{equation} \label{eq:e_SC}
\begin{split}
    a_{SC} &= \left(\frac{2}{r_{M1}} - \frac{V_{SC}^2}{GM_S} \right)^{-1}\\
    e_{SC} &= \sqrt{1 - \frac{h_{SC}^2} {a_{SC}  GM_S}}
\end{split}
\end{equation}

\noindent where $GM_S$ is the gravitational parameter of the Saturn and $a_{SC}$ is the semimajor axis of the SC orbit. Finally, we enforce the periapsis radius of the SC orbit to be smaller than the semimajor axis of the next moon $r_{M2}$ to find feasible transfers between moons.
\begin{equation} \label{eq:exit1}
\begin{split}
    r_{p,SC} &= a_{SC} (1-e_{SC}) \leq r_{M2}
\end{split}
\end{equation}


\noindent We then prune the feasible transfers to conform to the arrival $V_\infty$ bounds specified in Table \ref{tab:vinf_ub}. 

\textcolor{blue}{In this algorithm, the phasing of these inter-moon transfers is also considered. Before starting the DP, we precompute the set of feasible intermoon transfers by solving a Kepler's equation for the discretized range of $V_\infty$ and pump angle; this provides the rough estimates of the relationship between $\alpha_{dep}$ and the difference in transfer angle between the SC and the next moon in order to successfully match their position after the transfer, given the $V_\infty$. 
As we can compute the relative phase of the next moon from the SC at each node by carrying over the total ToF and initial phasing information, we can solve a nonlinear root finding problem with the initial guess based on the above estimates so that the SC and the next moon exactly matches their position. 
Using this method, we can compute the ToF and transfer angle of the adequate inter-moon transfer with minimum computation burden during the DP.} 

\subsection{Terminal condition of the path-finding: Enceladus orbit insertion}
At the end of the Saturn-moon tour (so-called \textit{endgame}), we include an Enceladus orbit insertion (EOI) $\Delta V$ to enter a $h=$ 100 km circular orbit. 

\begin{equation} \label{eq:EOI}
    \Delta V_{EOI} = \sqrt{V_{\infty}^2 + \frac{2GM_M}{r_M+h}} - \sqrt{\frac{GM_M}{r_M + h}}
\end{equation}

In this study, when $V_\infty$ of an uncompleted path becomes less than 450 m/s in the Enceladus tour, we calculate $\Delta V_{EOI}$ and store it in a memory for completed paths. These final solutions are filtered by 2D Pareto sorting based on the original two objectives (total ToF and total $\Delta V$) at the end of the optimization process. Note that the original uncompleted paths will still be revisited as the initial nodes of the next round of transfer until no reachable paths remain. 

\section{Optimization Result}

\begin{figure}[ht]
  \centering
  \includegraphics[width=0.65\linewidth]{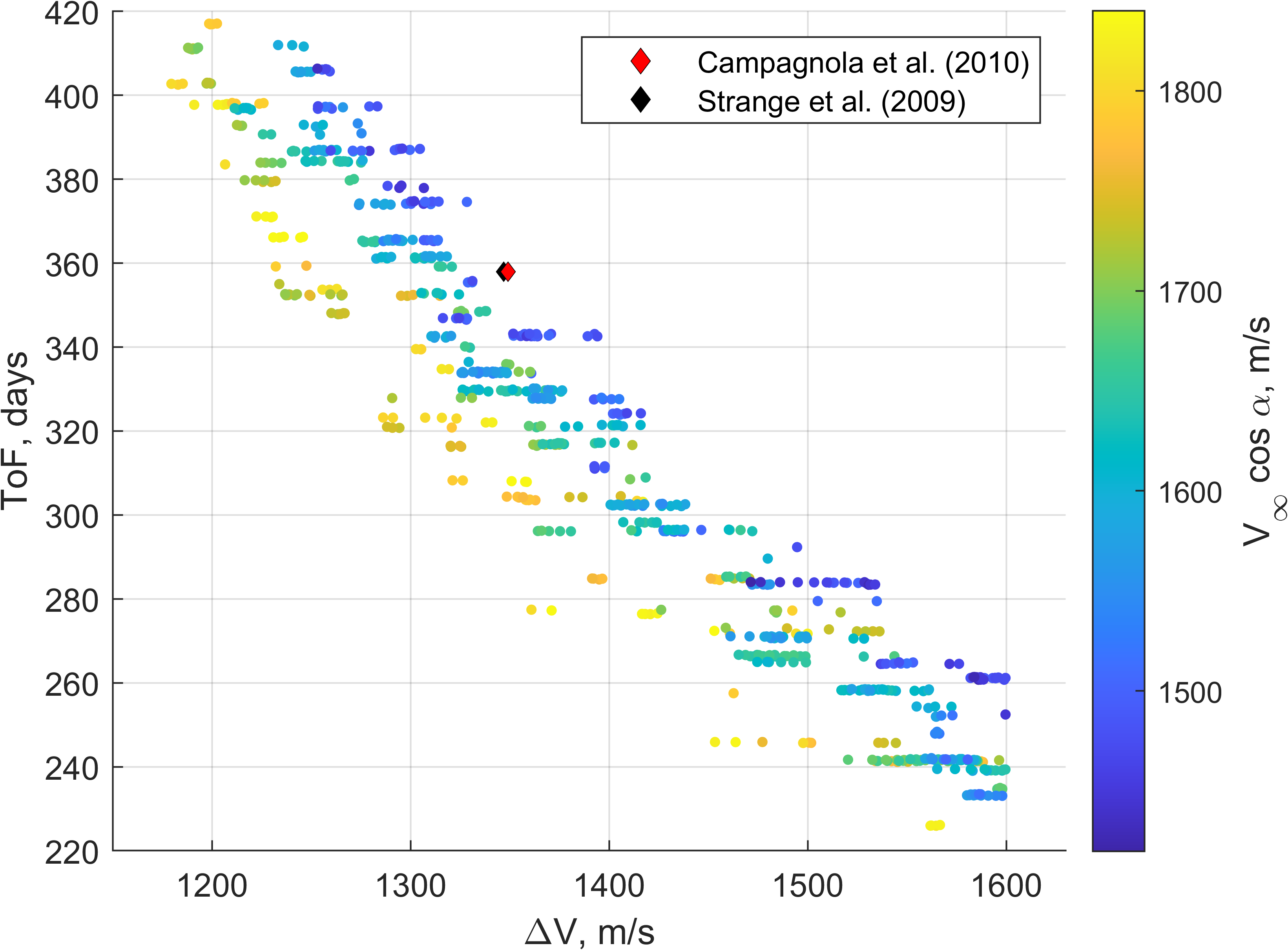}
  \caption{Initial states set after Titan tour. The color map indicates the transverse component of $V_\infty$ with respect to the next moon(i.e.,Rhea)'s velocity (i.e., $V_\infty \cos \alpha$).}
  \label{fig:titan_soln}
\end{figure}

\textcolor{blue}{The proposed method is applied to solve the tour design at Rhea, Dione, Tethys, and Enceladus. The initial states at Rhea (i.e., exiting states) at Titan is derived based on Star algorithm \cite{landau2022star}; it provides a single interplanetary trajectory from Earth to Saturn orbit insertion and numerous exiting states for the Titan tour, which are shown in Fig. \ref{fig:titan_soln}.} 
\textcolor{blue}{The red and black diamonds in the plot indicates solutions from past literature: Campagnola et al. \cite{campagnola2010enceladus}, and Strange et al. \cite{strange2009enceladus}, respectively. These values are provided by Star \cite{landau2022star}.}

The terminal state at Enceladus is set to a circular orbit of $h = 100$ km, and the total flight time is limited to 3 years. For the grid-based DP, we set an interval of the $V_\infty$ to $\Delta V_\infty = 30$ m/s. Also, for the VILT database, the interval of the $V_\infty$ is set to $\Delta V_\infty = 30$ m/s. The computation is performed on an i7-8550U @ 1.80GHz (4-core) CPU with 16 GB RAM, and the whole moon tour was completed in less than \textcolor{blue}{10 minutes}. 

\begin{figure}[H]
     \centering
     \begin{subfigure}[b]{0.49\textwidth}
         \includegraphics[width=\textwidth]{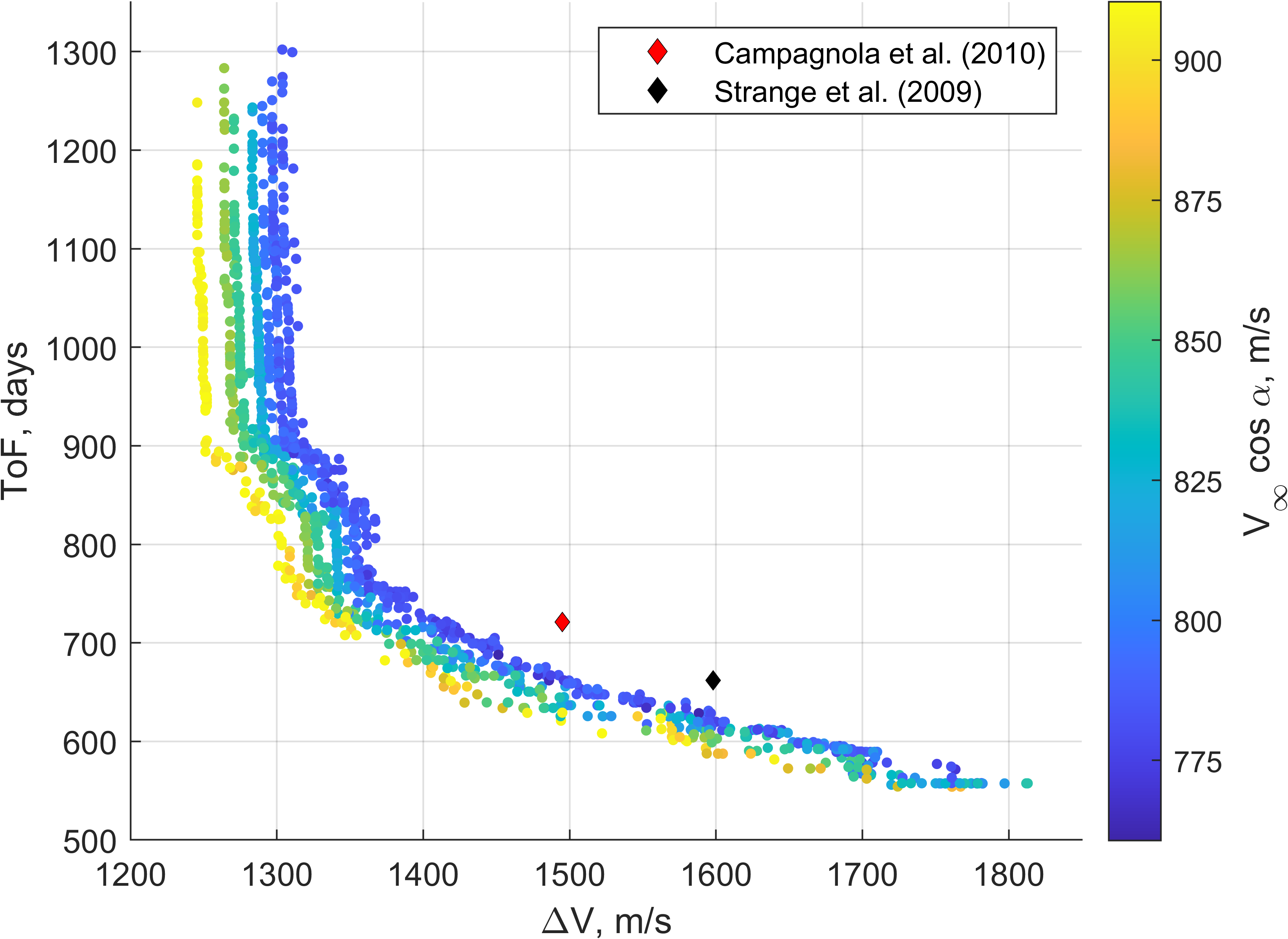}
     \end{subfigure}
     \begin{subfigure}[b]{0.49\textwidth}
         \includegraphics[width=\textwidth]{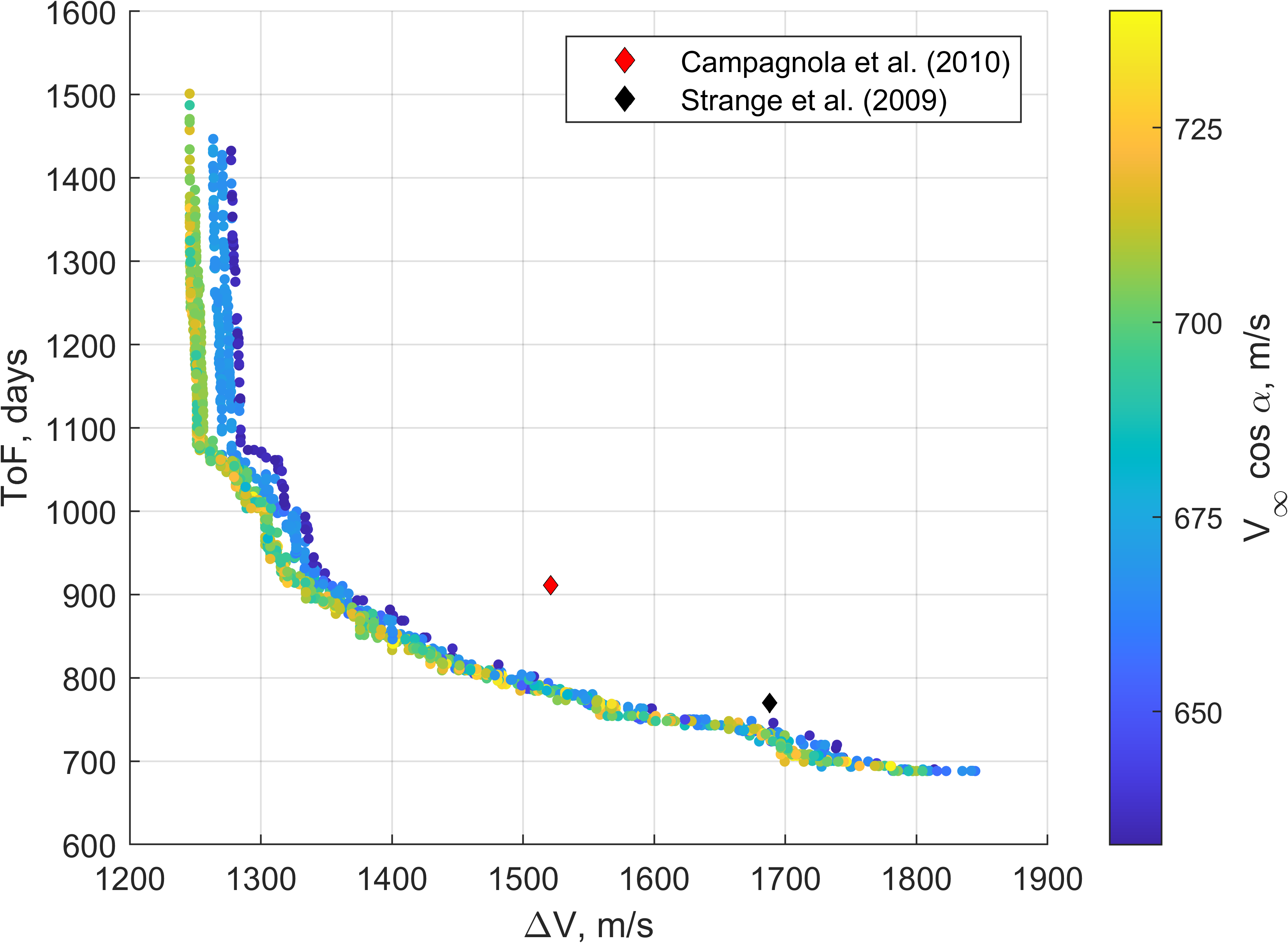}
     \end{subfigure}
     \newline
     \begin{subfigure}[b]{0.49\textwidth}
         \includegraphics[width=\textwidth]{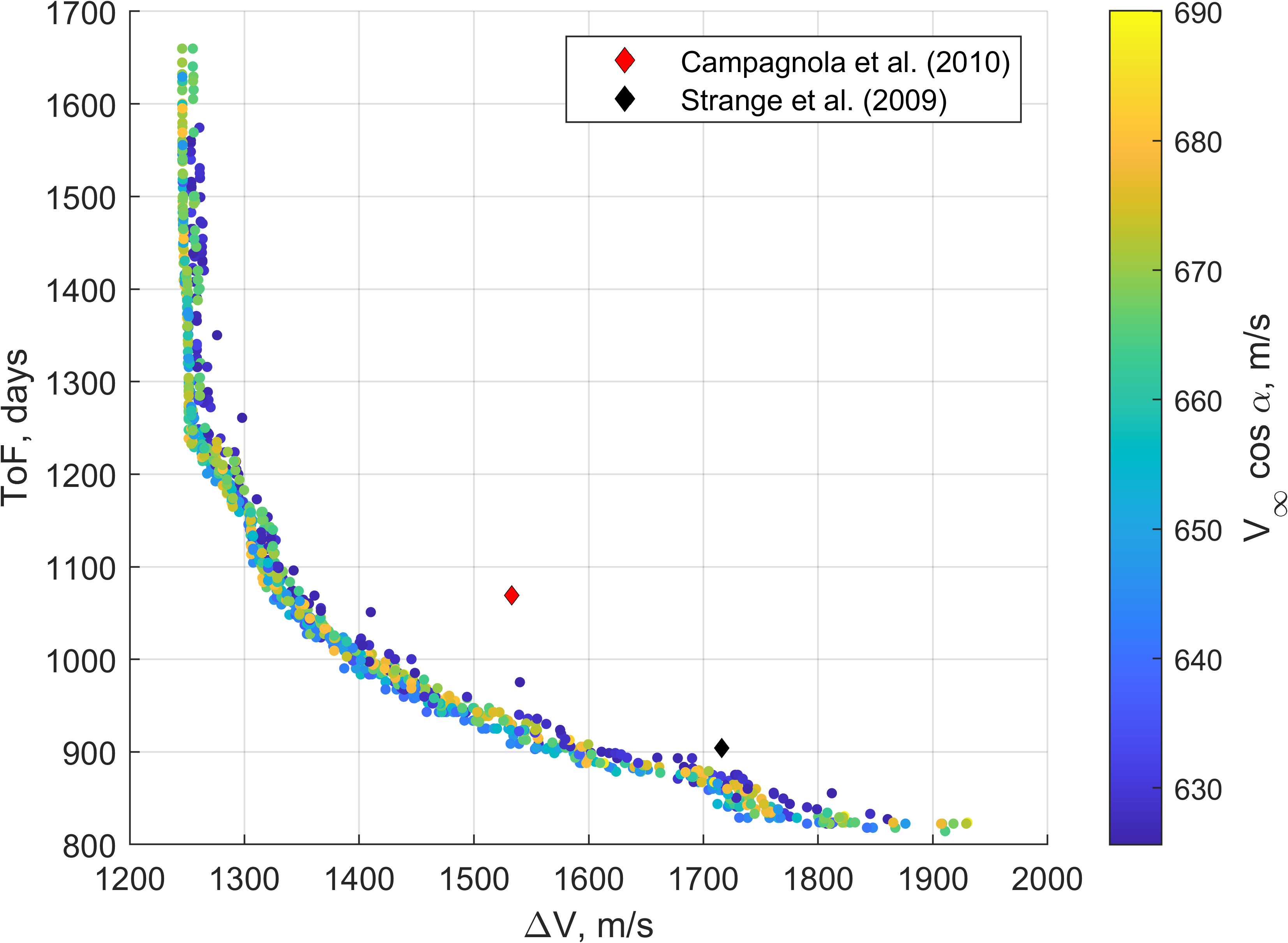}
     \end{subfigure}
     \begin{subfigure}[b]{0.49\textwidth}
         \includegraphics[width=\textwidth]{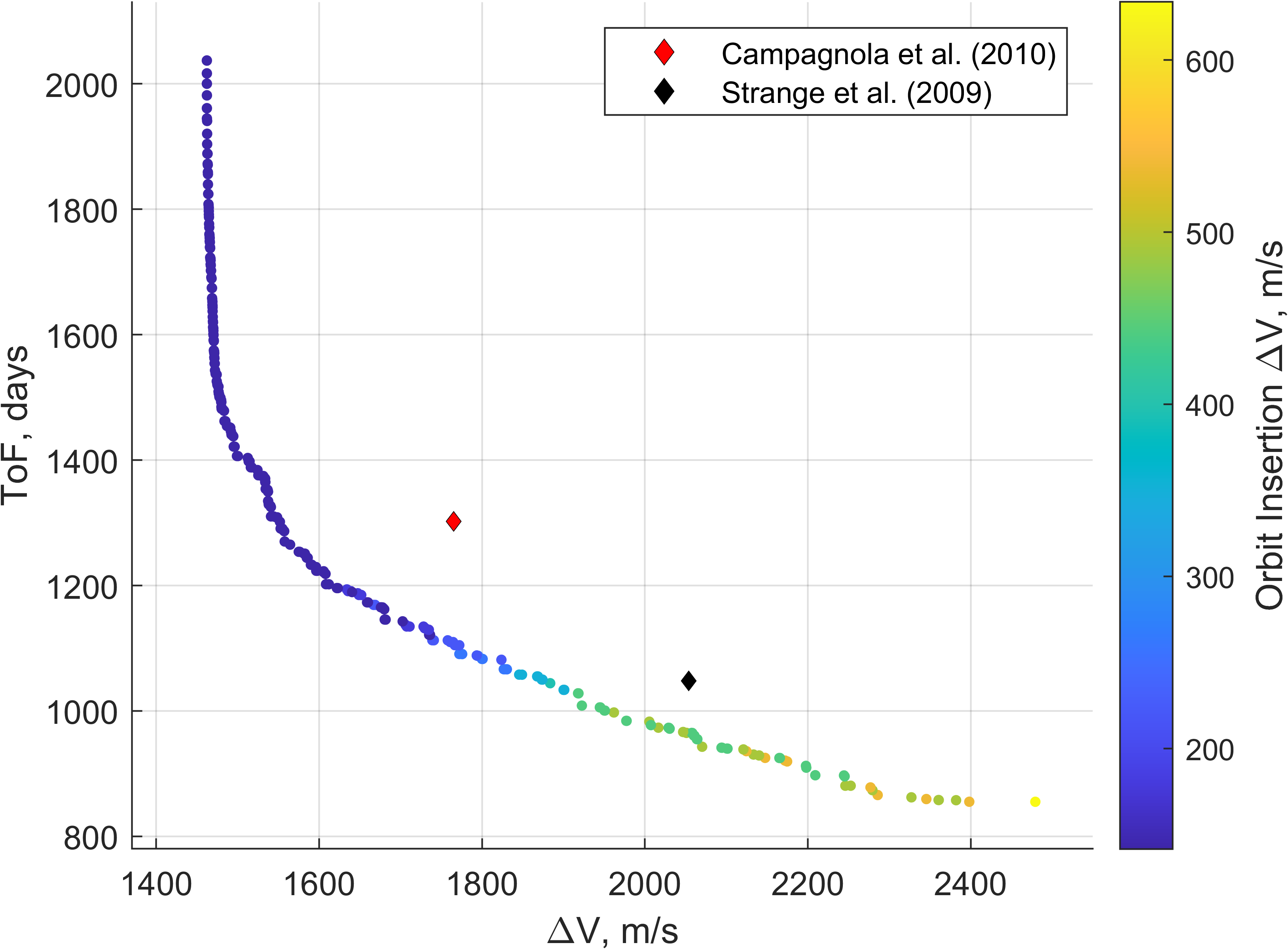}
     \end{subfigure}
     \caption{Feasible solution sets at the end of Rhea (top-left), Dione(top-right), and Tethys (bottom-left) tour, and 2D Pareto solution set at the end of Enceladus tour (bottom-right). The circled tour in the bottom right figure (Titan-Enceladus tour) is detailed in Figures \ref{fig:tour_pump_vinf} and \ref{fig:tour_traj}.}
\label{fig:tour_pareto}
\end{figure}

Figure \ref{fig:tour_pareto} represents the \textcolor{blue}{feasible solutions at the end of Rhea, Dione, and Tethys tour that are carried over to the next moon tour, and the final global Pareto solutions at Encedlaus moon tour. The color maps indicates the transverse component of the $V_\infty$ at the next moon for the first three moons, and the final $\Delta V$ for the EOI for the latter.}
To seek the trade space of the whole moon tour sequence from Titan to Enceladus, both $\Delta V$ ($x$-axis) and ToF ($y$-axis) are cumulative values from the \textcolor{blue}{SOI}. 
Therefore, each points in these plots include the information of not only the flyby sequence of that moon but all legs from Titan. 
Note that the final 2D Pareto front at Enceladus shown in the bottom-right panel of Figure \ref{fig:tour_pareto} includes the EOI cost, as shown in the color map. 
We can clearly confirm that the obtained solution sets at each moon tour improve the state of the art from the literature.
Our process automatically uncovers a collection of the trajectories that bound the previous solutions, with a ToF range of \textcolor{blue}{900 to 2000 days and total $\Delta V$ of 1400 to 2400 m/s.}

We choose one of the tours from the Pareto front (green circle in Figure \ref{fig:tour_pareto}) for validation in a full (ephemeris) model that includes ephemerides from JPL's Horizons system\footnote{\url{https://ssd-api.jpl.nasa.gov/doc/horizons.html} }, point-mass gravity from the five moons (Titan, Rhea, Dione, Tethys, and Enceladus), and Saturn zonal harmonics up to degree 10\footnote{Constants used to model the Saturnian system are available at \url{https://ssd.jpl.nasa.gov/ftp/eph/satellites/nio/LINUX_PC/sat441l.txt} }. 
First, the circular-coplanar solution (Figure \ref{fig:tour_pareto}) provides an initial guess for optimization in the \emph{Patched+} model, which uses Horizons ephemerides, pseudostate approximation for flybys, and secular Saturn-$J_2$ effects \cite{Kaela2019europa}. 
This initial guess is formulated by propagating the trajectories from each flyby to an apsidal maneuver (like Figure \ref{fig:tpbvp}), using the encounter dates, $V_\infty$, and pump angles of the circular-coplanar solution.
The approximate force model on the converged Patched+ (moon- and Saturn-centered) legs is replaced with the full model and optimized for minimum $\Delta V$ while maintaining the equivalent $V_\infty$ and pump angle necessary to connect to the previous and following Moons (or enter Enceladus orbit). 
The initial full-model result typically exceeds the estimated $\Delta V$ from circular-coplanar by several tens of m/s for each moon. We then apply the primer vector method \cite{lawden1963optimal} to add maneuvers that reduce the overall objective function to a local minimum. 
Optimization in Patched+ and full models occurs in a JPL in-house optimizer, ZoSo, that solves the Karush-Kuhn-Tucker conditions for optimality via a sequence of trust-region Newton steps with second-order Lagrange-multiplier update. ZoSo employs a collocation algorithm with adaptive mesh refinement to satisfy the SC dynamics. The Lagrange multipliers associated with the velocity collocation defects are equivalent to the primer vector\cite{landau2018efficient}, providing a convenient source to inform when additional maneuvers will improve the trajectory.

\begin{figure}[!ht]
     \centering
     \begin{subfigure}[b]{0.49\textwidth}
         \includegraphics[width=\textwidth]{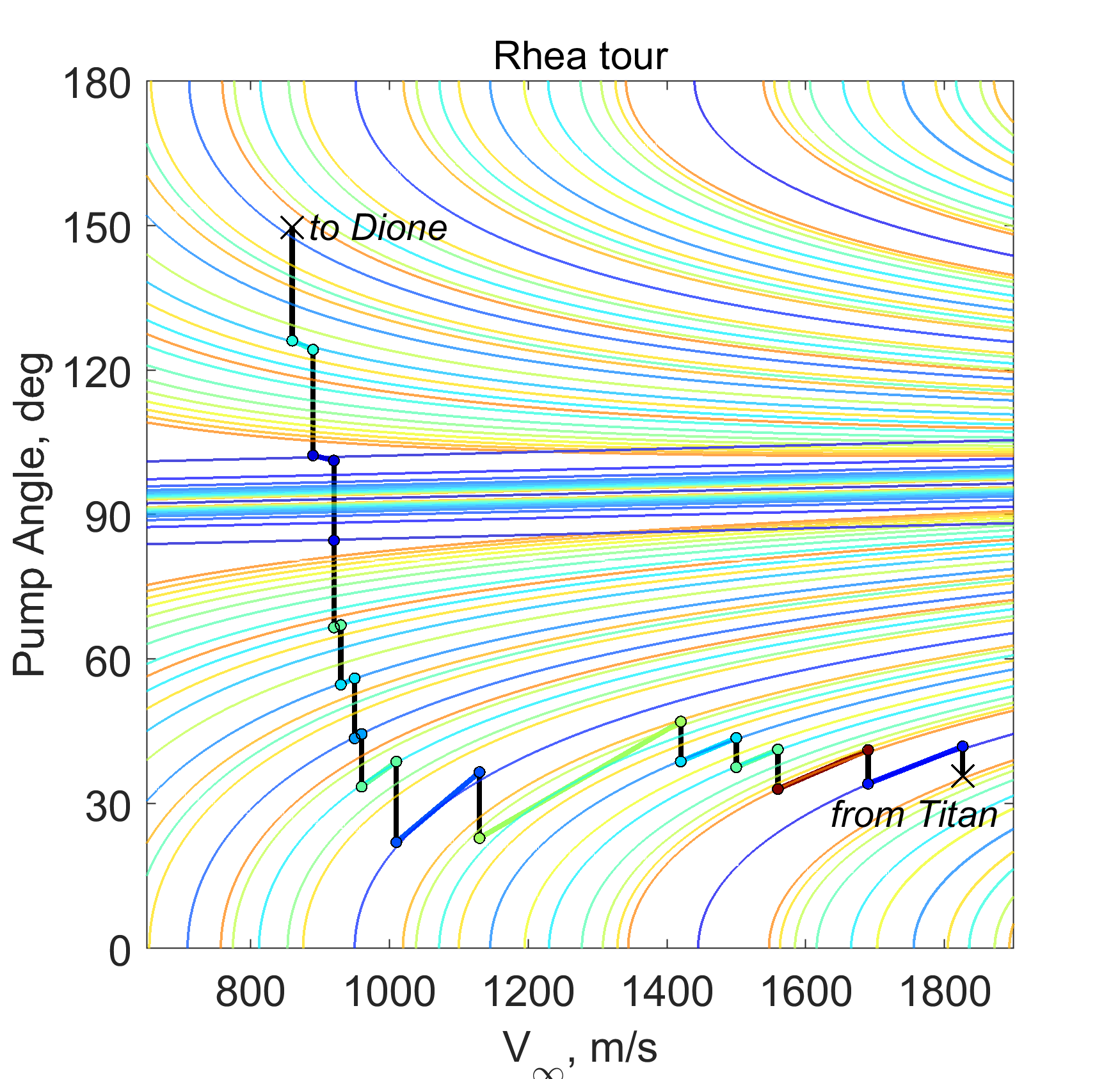}
     \end{subfigure}
     \begin{subfigure}[b]{0.49\textwidth}
         \includegraphics[width=\textwidth]{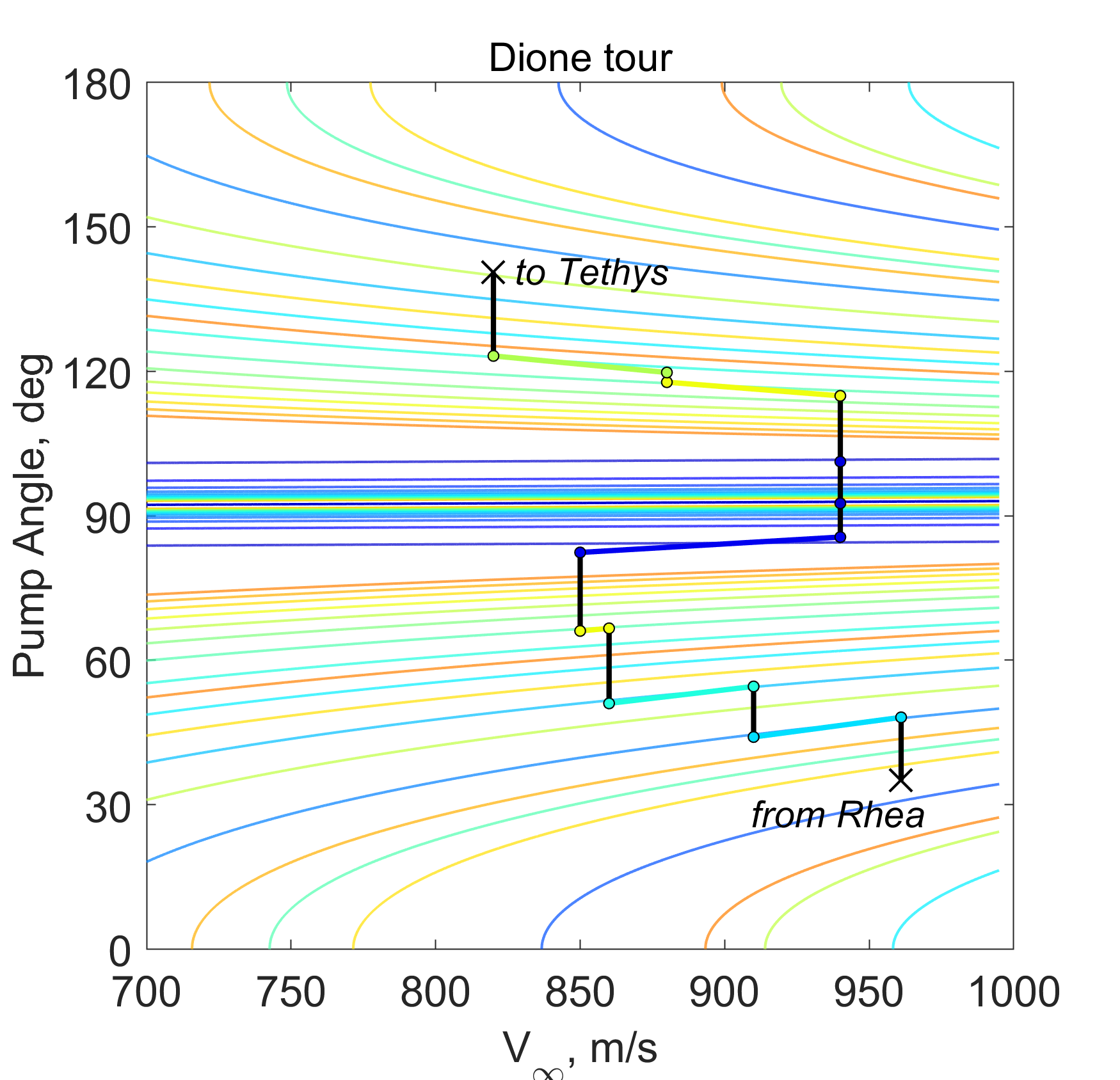}
     \end{subfigure}
     \newline
     \begin{subfigure}[b]{0.49\textwidth}
         \includegraphics[width=\textwidth]{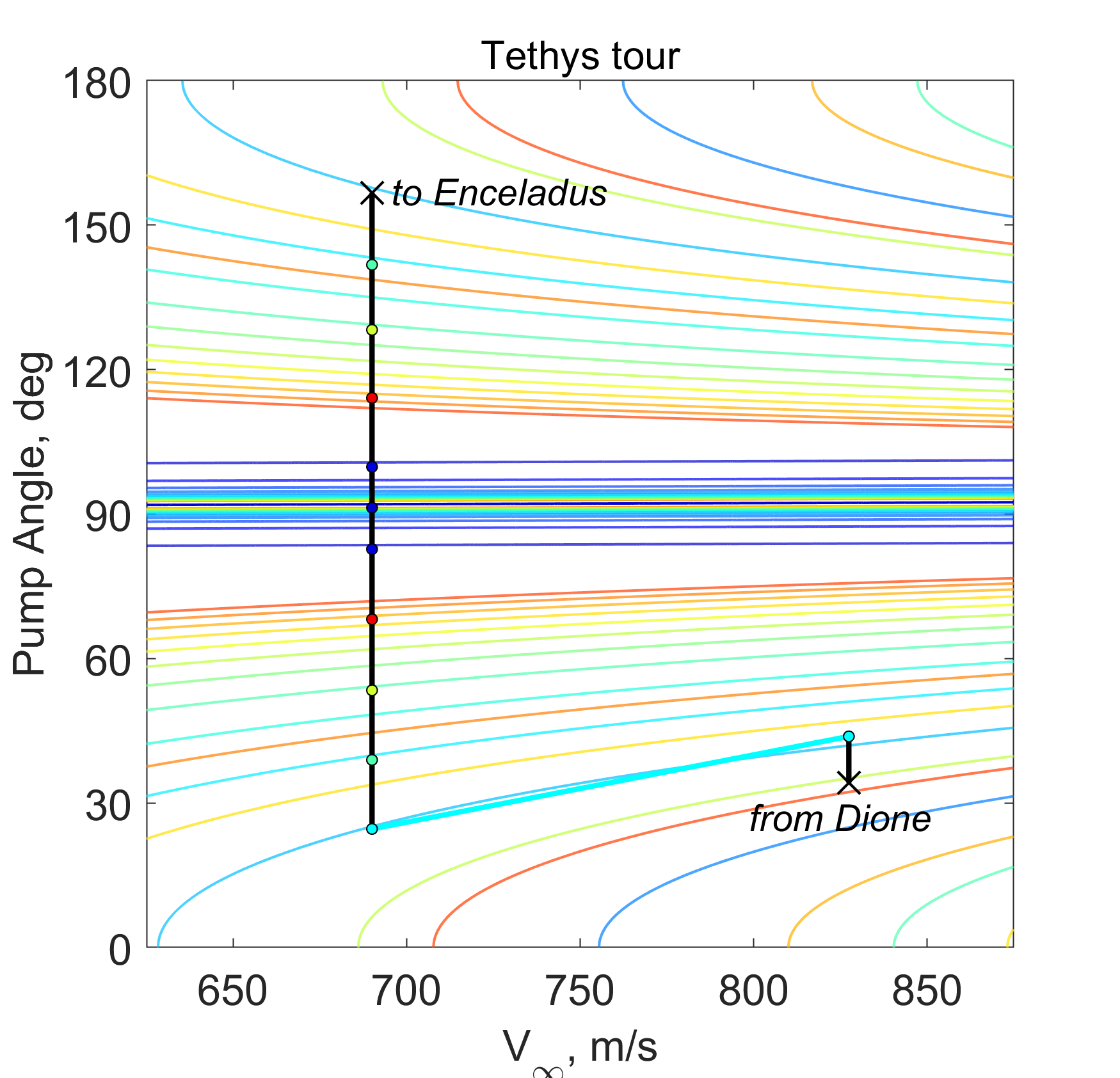}
     \end{subfigure}
     \begin{subfigure}[b]{0.49\textwidth}
         \includegraphics[width=\textwidth]{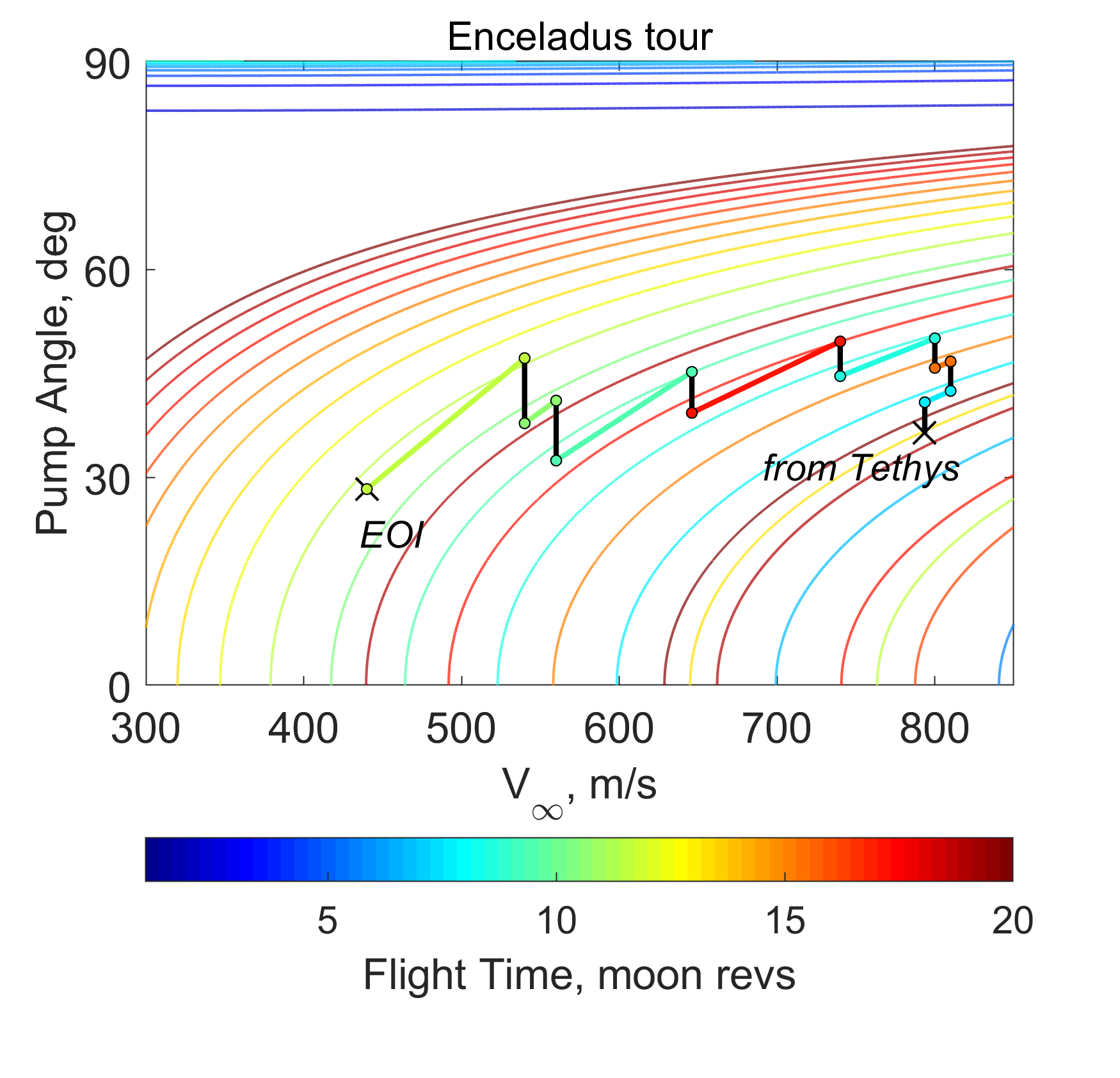}
     \end{subfigure}
     \caption{The combination of $V_\infty$ leveraging via impulsive $\Delta V$ (thick colored lines) and swings in pump via gravity assist (vertical black lines) travese the entire design space (circular-coplanar solution).}
\label{fig:tour_pump_vinf}
\end{figure}

\begin{figure}[!ht]
  \centering
  \includegraphics[width=0.8\linewidth]{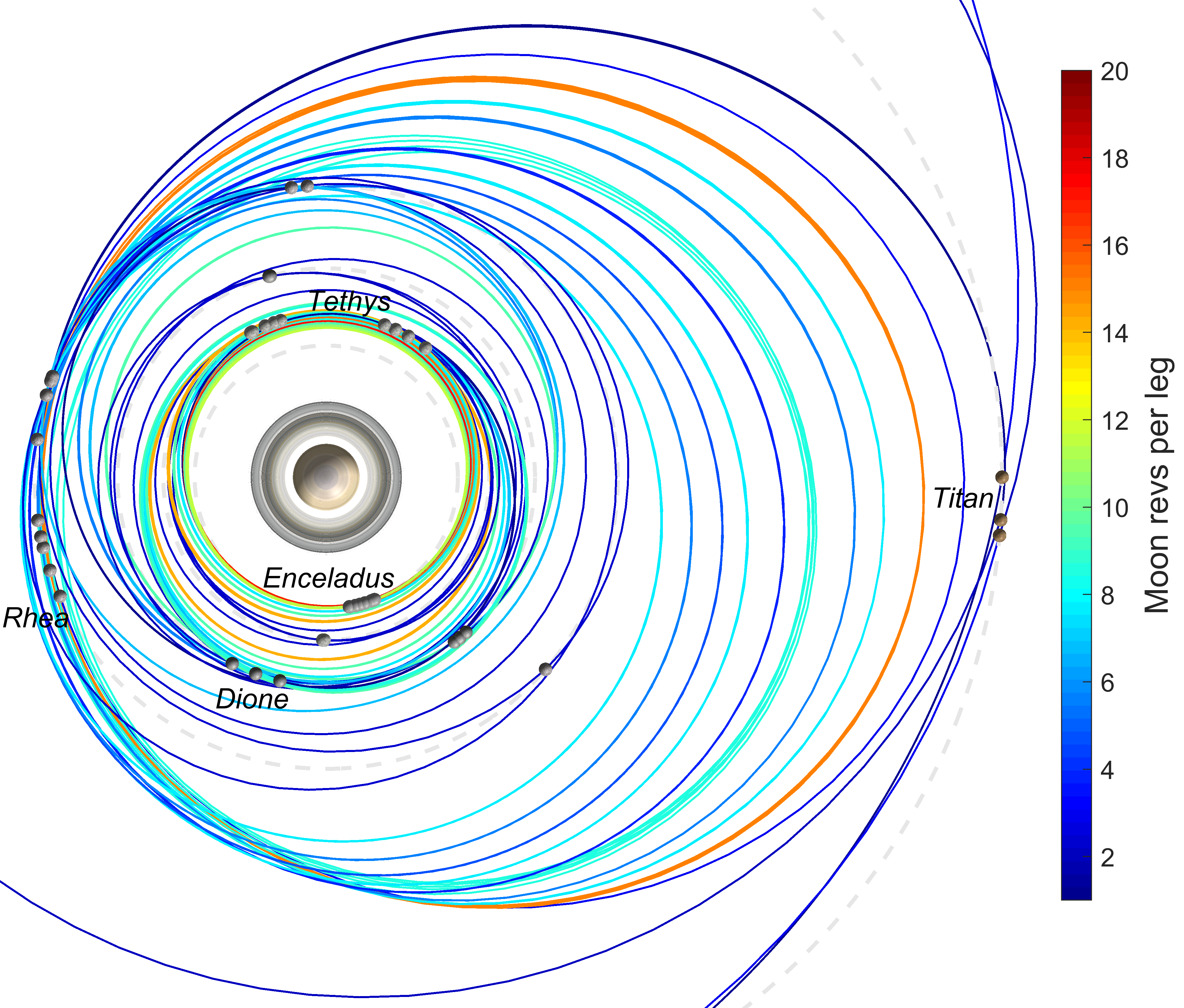}
  \caption{End-to-end tour from a 2:1 resonance with Titan to Enceladus orbit insertion (circular-coplanar solution).}
  \label{fig:tour_traj}
\end{figure}

The tour of our representative circular-coplanar solution from Figure \ref{fig:tour_pareto} is shown in Figures \ref{fig:tour_pump_vinf} and \ref{fig:tour_traj}. Tables \ref{tab:Tit_tour}--\ref{tab:Enc_tour} contain additional details of the individual flybys and transfers for comparison between the circular-coplanar initial estimate and the converged full-model solution. The initial and final flybys of each moon are constrained to maintain continuity with the rest of the tour, while the remaining flybys shift to minimize $\Delta V$ while maintaining altitude above 50 km. The input sequence of transfer types locks each moon tour into a local minimum. The summary contained in Table \ref{tab:sum_tour} demonstrates that the addition of optimal maneuvers in the full model further reduces total $\Delta V$ by 90 m/s while maintaining the total flight time. This refined solution makes a strong case that our Pareto solution sets a reliable benchmark for design of sequential moon tours at Saturn.

\begin{table}[!ht]
\centering
\begin{tabular}{c c|c c c c|c c c c} 
\hline
Flyby & Transfer & \multicolumn{4}{c}{Circular-Coplanar} & \multicolumn{4}{c}{Full Ephemeris}\\
Titan & Type & ToF & Alt. & $V_\infty$ & $\Delta V$ & ToF & Alt. & $V_\infty$ & $\Delta V$ \\ 
\# & $M:N^{p,q}$ & (day) & (km) & (m/s) & (m/s) & (day) & (km) & (m/s) & (m/s) \\ 
\hline\hline
1 & $2:1^{+,+}$ & 31.67 & 9080 & 1460 & 16.2 & 31.34 & 1600 & 1460 & 20.3 \\
2 & $1:1^{+,+}$ & 16.00 & 4884 & 1350 & 14.0 & 15.72 & 3497 & 1311 & 0.0 \\
3 & Rhea & --- & 1600 & 1320 & --- & --- & 1600 & 1320 & --- \\
\hline
\end{tabular}
\caption{Titan tour.}
\label{tab:Tit_tour}
\end{table}
\begin{table}[!ht]
\centering
\begin{tabular}{c c|c c c c|c c c c} 
\hline
Flyby & Transfer & \multicolumn{4}{c}{Circular-Coplanar} & \multicolumn{4}{c}{Full Ephemeris}\\
Rhea & Type & ToF & Alt. & $V_\infty$ & $\Delta V$ & ToF & Alt. & $V_\infty$ & $\Delta V$ \\ 
\# & $M:N^{p,q}$ & (day) & (km) & (m/s) & (m/s) & (day) & (km) & (m/s) & (m/s) \\ 
\hline\hline
1 & $2:1^{+,+}$ & 8.97 & 55 & 1826 & 16.9 & 9.03 & 50 & 1826 & 0.0 \\
2 & $15:8^{+,+}$ & 67.71 & 65 & 1690 & 16.6 & 67.68 & 158 & 1821 & 25.4 \\
3 & $7:4^{+,+}$ & 31.60 & 58 & 1560 & 8.1 & 31.53 & 50 & 1567 & 0.0 \\
4 & $5:3^{+,+}$ & 22.55 & 461 & 1500 & 11.3 & 22.38 & 194 & 1428 & 17.4 \\
5 & $8:5^{+,+}$ & 35.95 & 219 & 1420 & 39.3 & 36.14 & 50 & 1190 & 9.7 \\
6 & $3:2^{+,+}$ & 13.44 & 127 & 1130 & 16.0 & 13.48 & 50 & 1181 & 11.1 \\
7 & $7:5^{+,+}$ & 31.59 & 123 & 1010 & 7.3 & 31.36 & 50 & 1104 & 36.1 \\
8 & $4:3^{+,+}$ & 18.07 & 826 & 960 & 1.6 & 18.05 & 158 & 861 & 0.0 \\
9 & $5:4^{+,+}$ & 22.58 & 639 & 950 & 3.9 & 22.58 & 435 & 853 & 0.0 \\
10 & $7:6^{+,-}$ & 31.63 & 713 & 930 & 2.5 & 30.39 & 532 & 858 & 0.0 \\
11 & $1:1^{-,+}$ & 6.51 & 213 & 920 & 0.0 & 6.49 & 352 & 900 & 0.0 \\
12 & $1:1^{+,-}$ & 6.21 & 316 & 920 & 11.7 & 6.19 & 285 & 859 & 0.0 \\
13 & $6:7^{-,+}$ & 27.15 & 63 & 890 & 7.0 & 26.21 & 137 & 896 & 0.0 \\
14 & Dione & --- & 50 & 860 & --- & --- & 50 & 860 & --- \\
\hline
\end{tabular}
\caption{Rhea tour.}
\label{tab:Rhea_tour}
\end{table}
\begin{table}[!ht]
\centering
\begin{tabular}{c c|c c c c|c c c c} 
\hline
Flyby & Transfer & \multicolumn{4}{c}{Circular-Coplanar} & \multicolumn{4}{c}{Full Ephemeris}\\
Dione & Type & ToF & Alt. & $V_\infty$ & $\Delta V$ & ToF & Alt. & $V_\infty$ & $\Delta V$ \\ 
\# & $M:N^{p,q}$ & (day) & (km) & (m/s) & (m/s) & (day) & (km) & (m/s) & (m/s) \\ 
\hline\hline
1 & $5:4^{+,+}$ & 13.67 & 52 & 961 & 8.7 & 13.62 & 50 & 961 & 14.7 \\
2 & $6:5^{+,+}$ & 16.41 & 316 & 910 & 9.6 & 16.42 & 184 & 839 & 10.4 \\
3 & $9:8^{+,-}$ & 24.64 & 69 & 860 & 2.5 & 23.89 & 50 & 824 & 0.0 \\
4 & $1:1^{-,+}$ & 3.92 & 51 & 850 & 36.0 & 3.91 & 60 & 857 & 0.0 \\
5 & $1:1^{+,+}$ & 2.74 & 708 & 940 & 0.0 & 2.73 & 628 & 802 & 0.0 \\
6 & $1:1^{+,-}$ & 3.78 & 453 & 940 & 0.0 & 3.77 & 588 & 803 & 0.0 \\
7 & $9:10^{-,+}$ & 24.69 & 51 & 940 & 16.7 & 24.01 & 103 & 856 & 0.0 \\
8 & $8:9^{+,+}$ & 21.95 & 4731 & 880 & 15.1 & 21.89 & 1645 & 773 & 0.0 \\
9 & Tethys & --- & 50 & 820 & --- & --- & 50 & 820 & --- \\
\hline
\end{tabular}
\caption{Dione tour.}
\label{tab:Dio_tour}
\end{table}
\begin{table}[!ht]
\centering
\begin{tabular}{c c|c c c c|c c c c} 
\hline
Flyby & Transfer & \multicolumn{4}{c}{Circular-Coplanar} & \multicolumn{4}{c}{Full Ephemeris}\\
Tethys & Type & ToF & Alt. & $V_\infty$ & $\Delta V$ & ToF & Alt. & $V_\infty$ & $\Delta V$ \\ 
\# & $M:N^{p,q}$ & (day) & (km) & (m/s) & (m/s) & (day) & (km) & (m/s) & (m/s) \\ 
\hline\hline
1 & $6:5^{+,+}$ & 11.27 & 124 & 827 & 21.1 & 11.33 & 50 & 827 & 0.0 \\
2 & $7:6^{+,+}$ & 13.23 & 75 & 690 & 0.0 & 13.22 & 64 & 828 & 10.6 \\
3 & $9:8^{+,+}$ & 17.01 & 68 & 690 & 0.0 & 17.01 & 50 & 819 & 27.0 \\
4 & $14:13^{+,-}$ & 26.46 & 55 & 690 & 0.0 & 25.83 & 50 & 792 & 0.0 \\
5 & $1:1^{-,+}$ & 2.68 & 64 & 690 & 0.0 & 2.70 & 130 & 842 & 0.0 \\
6 & $1:1^{+,+}$ & 1.89 & 534 & 690 & 0.0 & 1.88 & 176 & 785 & 0.0 \\
7 & $1:1^{+,-}$ & 2.64 & 550 & 690 & 0.0 & 2.61 & 155 & 786 & 0.0 \\
8 & $14:15^{-,+}$ & 26.47 & 75 & 690 & 0.0 & 25.93 & 171 & 844 & 20.1 \\
9 & $9:10^{+,+}$ & 17.02 & 84 & 690 & 0.0 & 17.01 & 50 & 723 & 0.0 \\
10 & $7:8^{+,+}$ & 13.24 & 116 & 690 & 0.0 & 13.24 & 50 & 722 & 6.6 \\
11 & Enceladus & --- & 50 & 690 & --- & --- & 50 & 690 & --- \\
\hline
\end{tabular}
\caption{Tethys tour.}
\label{tab:Tet_tour}
\end{table}
\begin{table}[!ht]
\centering
\begin{tabular}{c c|c c c c|c c c c} 
\hline
Flyby & Transfer & \multicolumn{4}{c}{Circular-Coplanar} & \multicolumn{4}{c}{Full Ephemeris}\\
Enceladus & Type & ToF & Alt. & $V_\infty$ & $\Delta V$ & ToF & Alt. & $V_\infty$ & $\Delta V$ \\ 
\# & $M:N^{p,q}$ & (day) & (km) & (m/s) & (m/s) & (day) & (km) & (m/s) & (m/s) \\ 
\hline\hline
1 & $7:6^{+,+}$ & 9.61 & 32 & 794 & 2.8 & 9.58 & 50 & 794 & 10.9 \\
2 & $15:13^{+,+}$ & 20.58 & 34 & 810 & 1.8 & 20.57 & 50 & 730 & 11.9 \\
3 & $8:7^{+,+}$ & 10.96 & 40 & 800 & 10.9 & 10.97 & 50 & 689 & 0.0 \\
4 & $17:15^{+,+}$ & 23.30 & 37 & 740 & 16.5 & 23.31 & 50 & 682 & 8.2 \\
5 & $9:8^{+,+}$ & 12.32 & 66 & 646 & 14.2 & 12.31 & 50 & 644 & 11.5 \\
6 & $10:9^{+,+}$ & 13.71 & 30 & 560 & 3.3 & 13.72 & 50 & 570 & 21.3 \\
7 & $11:10^{+,+}$ & 15.05 & 25 & 540 & 16.5 & 14.97 & 50 & 542 & 30.6 \\
8 & EOI & --- & 100 & 440 & 341 & --- & 100 & 386 & 293 \\
\hline
\end{tabular}
\caption{Enceladus tour.}
\label{tab:Enc_tour}
\end{table}
\begin{table}[!ht]
\centering
\begin{tabular}{c|c c|c c|c c} 
\hline
Moon & \multicolumn{2}{c}{Circular-Coplanar} & \multicolumn{2}{c}{Patched+} & \multicolumn{2}{c}{Full Ephemeris}\\
Tour & ToF (day) & $\Delta V$ (m/s) & ToF (day) & $\Delta V$ (m/s) & ToF (day) & $\Delta V$ (m/s) \\ 
\hline\hline
Titan & 47.7 & 30.2 &  47.7 & 20.0 & 47.1 & 20.3  \\
Rhea & 324.0 & 142.1 & 322.0 & 126.1 & 321.5 & 99.7  \\
Dione & 111.8 & 88.6 & 111.3 & 26.3 & 110.2 & 25.2  \\
Tethys & 131.9 & 21.1 & 130.8 & 61.1 & 130.8 & 64.3  \\
Enceladus & 105.5 & 66.0 & 105.4 & 93.6 & 105.4 & 94.3  \\
EOI & --- & 341 & --- & 282 & --- & 293  \\
\hline
Total &  721 & 689 & 717 & 609 & 715 & 597  \\
\hline
\end{tabular}
\caption{Summary comparison of tours and design models.}
\label{tab:sum_tour}
\end{table}

\section{Conclusion}
Multi-objective optimization of multi-moon tours in the Saturnian system can be automated to produce a trade space of total flight time versus $\Delta V$ to reach the Enceladus orbit. Several approximation and discretization techniques simplify the computational burden of completing a global search within a few hours. These simplified results are sufficiently accurate to provide the initial guess for  high-fidelity trajectory optimization. New insights to the available trade space to reach Enceladus expand the possibilities to search for life outside of Earth.  

\section{Acknowledgement}
Yuji Takubo thanks Keidai Iiyama for providing a code for the Tisserand graphs. The initial research was funded by JPL's Research and Technology Development program and investigated by Damon Landau, Stefano Campagnola, Reza Karimi, and Nathan Strange. A portion of the research described in this paper was carried about at the Jet Propulsion Laboratory, California Institute of Technology, under a contract with the National Aeronautics and Space Administration.

\appendix
\section*{Appendix}
\subsection{Initial guesses of the leg optimization}
\textit{Case 1: $p=q$}.
We can derive an analytical solution for this case. 
The SC and the moon in this case make an exact integer revolutions, as the departure and the arrival happen at the same position; $t_f$ is exactly the time for $N$ revolutions of Saturn (i.e., total transfer angle is $\theta \equiv \theta_1 + \theta_2 = 2 \pi N$), so we can compute it from Table \ref{tab:moon_param}. 
From the Kepler's third law and vis-viva equation, we can obtain the SC velocity when it encounters the moon. 
\begin{equation} \label{eq:v_SC_ballistic}
\begin{split}
    V_{SC} &= \sqrt{GM_S \left(\frac{2}{r_M} - \frac{1}{a_{SC}}\right)} \\
           &= \sqrt{GM_S \left(\frac{2}{r_M} - \frac{1}{r_M} \left(\frac{N}{M}\right)^{\frac{2}{3}}\right)} 
\end{split}
\end{equation}

\noindent Using the law of cosine, we can derive the pump angle as  
\begin{equation} \label{eq:alpha}
\begin{split}
    \alpha = \cos^{-1} \left( \dfrac{V_{SC}^2 - V_M^2 - V_\infty^2}{2 V_M V_\infty} \right)
\end{split}
\end{equation}

\noindent The angular momentum and eccentricity are then obtained via Eqs. \ref{eq:h_sc} and \ref{eq:e_SC}, providing the absolute value of the true anomaly where the SC encounters the moon. 

%
%
\begin{equation} \label{eq:f_SC}
\begin{split}
    f(\alpha) &= \cos^{-1} \left( \frac{1}{e(\alpha)} \left(\frac{h_{SC}(\alpha)^2}{r_M GM_S} - 1 \right) \right), \quad 0 \le f(\alpha) < \pi
\end{split}
\end{equation}

\textit{Case 2: $p \neq q$}. 
There is no analytical solution when the leg geometry is symmetric, and we develop a numerical solution method. This is solved via a root finding problem for the pump angle by constraining the flight time of the SC and the moon to be equal.
\begin{equation} \label{eq:root_finding}
\begin{split}
    t_{f, SC}(\alpha) - t_{f,M}(\alpha) = 0
\end{split}
\end{equation}

\noindent We obtain the ToF of the SC from $\alpha$ via Eqs. \ref{eq:V_sc}, \ref{eq:h_sc}, \ref{eq:e_SC}, \ref{eq:f_SC}, and the eccentric anomaly. 
\begin{equation} \label{eq:t_f_SC}
\begin{split}
    E &= 2 \tan^{-1} \sqrt{\left(\frac{1-e}{1+e}\right) \tan \left(\frac{f}{2}\right)}\\
    t_{f,SC}(\alpha) &= \sqrt{a_{SC}^3/GM_S} \left[ 2 \pi (N+\Delta) + (E(f_p) - e \sin E(f_p)) - (E(f_q) - e \sin E(f_q)) \right]
\end{split}
\end{equation}

\noindent Similarly, the total transfer angle for the SC is
\begin{equation} \label{eq:total_trans_angle}
    \theta \equiv \theta_1 + \theta_2 = 2 \pi (N+\Delta)  + f_q (\alpha) - f_p(\alpha)
\end{equation}

\noindent where $f_p = pf$ and $f_q = qf$ indicate the true anomaly at departure and arrival. The ToF of the moon is determined from velocity on a circular orbit
\begin{equation} \label{eq:t_f_M}
\begin{split}
    t_{f,M}(\alpha) = \left(2 \pi (M+\Delta) + f_q (\alpha) - f_p(\alpha)\right)r_M/V_M
\end{split}
\end{equation}

\noindent where $\Delta$ corrects for the number of SC periapses on the different resonance families.
\begin{equation*}
\Delta = 
\begin{cases}
    0 & \text{if } (M > N) \text{ or } (M=N \text{ and } p=-1, q=+1)\\
    p & \text{otherwise} \\
\end{cases}
\end{equation*}


\noindent With Eq. \ref{eq:t_f_SC} and \ref{eq:t_f_M}, Eq. \ref{eq:root_finding} is solved for $\alpha$ via Newton's method. To take the numerical derivative, we use the complex step differentiation \cite{Martins2003complex}.  
For the initial guess of $\alpha$, the solution of Case 1 that has the same resonance ratio $M:N$ is adopted. 

Finally, for both Case 1 and Case 2, the conversion from $f(\alpha)$ to $\theta_1$ and $\theta_2$ is
\begin{equation} \label{eq:f_to_thetas}
\begin{split}
    \theta_1 &= \pi \left(1+\Delta+2\cdot\text{floor}(N/2)\right) - f_p(\alpha)\\
    \theta_2 &= \theta-\theta_1    
\end{split}
\end{equation}

\bibliographystyle{AAS_publication}   
\bibliography{references}   

\end{document}